\documentclass[12pt]{article}
\usepackage{amsmath,amsthm,amssymb,times,titlesec}
\usepackage{graphicx}
\usepackage{pdfsync}
\usepackage{float}
\numberwithin{equation}{section}

\textheight 
            650pt
\textwidth 
           460pt

\titleformat{\subsubsection}
{\normalfont\fontsize{13}{14}\selectfont\itshape}{\thesubsubsection}{1em}{}

\newcommand\scaleddot{\scalebox{.89}{.}}
\makeatletter
\renewcommand{\ddddot}[1]{%
  {\mathop{\kern\z@#1}\limits^{\makebox[0pt][c]{\vbox to-2\ex@{\kern-\tw@\ex@\hbox{\normalfont\scaleddot\kern-0.5pt\scaleddot\kern-0.5pt\scaleddot\kern-0.5pt\scaleddot}\vss}}}}}

\newcommand{\nn}{|{\mskip-2mu}|{\mskip-2mu}|}

\newcommand{\C}{\mathbb{C}}
\newcommand{\R}{\mathbb{R}}

\newcommand{\N}{\mathbb{N}}

\newcommand{\ee}{\mathrm{e}}
\newcommand{\ii}{\mathrm{i}}

\newcommand{\dk}{\, \mathrm{d}k}

\newcommand{\dx}{\, \mathrm{d}x}

\newcommand{\dtildey}{\, \mathrm{d}\tilde{y}}
\newcommand{\dz}{\, \mathrm{d}z}

\newcommand{\FF}{{\mathcal F}}
\newcommand{\GG}{{\mathcal G}}

\newcommand{\JJ}{{\mathcal J}}
\newcommand{\KK}{{\mathcal K}}
\newcommand{\LL}{{\mathcal L}}
\newcommand{\KKp}{{\mathcal K}^\prime}
\newcommand{\LLp}{{\mathcal L}^\prime}

\newcommand{\NN}{{\mathcal N}}

\renewcommand{\SS}{{\mathcal S}}

\newcommand{\XX}{{\mathcal X}}

\newcommand{\ZZ}{{\mathcal Z}}

\DeclareMathOperator{\re}{Re}
\DeclareMathOperator{\im}{Im}
\DeclareMathOperator{\sgn}{sgn}

\newtheorem{theorem}{Theorem}[section]
\newtheorem{lemma}[theorem]{Lemma}
\newtheorem{proposition}[theorem]{Proposition}
\newtheorem{corollary}[theorem]{Corollary}
\newtheorem{remark}[theorem]{Remark}
\newtheorem{definition}[theorem]{Definition}
\theoremstyle{definition}

\renewcommand{\i}{\mathrm{i}}

\newcounter{count}

\title{Fully localised three-dimensional gravity-capillary solitary waves on water of infinite depth}

\author{B. Buffoni\thanks{Institut de math\'{e}matiques, Station 8, \'{E}cole polytechnique f\'{e}d\'{e}rale, 1015 Lausanne, Switzerland
}\and
M. D. Groves\footnote{Fachrichtung Mathematik, Universit\"{a}t des Saarlandes,
Postfach 151150, 66041 Saarbr\"{u}cken, Germany}
\and
E. Wahl\'{e}n\footnote{Centre for Mathematical Sciences, Lund University, PO Box 118, 22100 Lund, Sweden}}
\date{}

\begin{document}

\maketitle

\begin{abstract}
\emph{Fully localised solitary waves} are travelling-wave solutions
of the three-dimensional gravity-capillary water wave problem which
decay to zero in every horizontal spatial direction. Their existence for
water of finite depth has recently been established, and in this article
we present an existence theory for water of infinite depth. The governing
equations are reduced to a perturbation of the two-dimensional
nonlinear Schr\"{o}dinger equation, which admits a family
of localised solutions. Two of these solutions are symmetric in both
horizontal directions and an application of a suitable version of the
implicit-function theorem shows that they persist under perturbations.
\end{abstract}

\section{Introduction}

Three-dimensional gravity-capillary water waves on the surface of
a body of water of infinite depth are described by the Euler equations
in a domain bounded above by a free
surface $\{y=\eta(x,z,t)\}$, where the
function $\eta$ depends upon the two horizontal spatial directions
$x$, $z$ and time $t$.
In terms of an Eulerian velocity potential $\varphi$ and in dimensionless coordinates,
the mathematical problem is to solve Laplace's equation
\begin{equation}
\varphi_{xx} + \varphi_{yy} + \varphi_{zz} = 0, \qquad\qquad -\infty<y<\eta, \label{SWW 1}
\end{equation}
with boundary conditions
\begin{eqnarray}
\varphi_{y} & \rightarrow & \parbox{50mm}{$0$,} y \rightarrow -\infty, \label{SWW 2} \\
\eta_t & = & \parbox{50mm}{$\varphi_{y} - \eta_x\varphi_x - \eta_z\varphi_z$,}  y=\eta, \label{SWW 3}
\end{eqnarray}
and
\begin{equation}
\varphi_t  =  -\frac{1}{2}(\varphi_x^2+\varphi_{y}^2+\varphi_z^2) -\eta
+ \left[\frac{\eta_x}{\sqrt{1+\eta_x^2+\eta_z^2}}\right]_x
+ \left[\frac{\eta_z}{\sqrt{1+\eta_x^2+\eta_z^2}}\right]_z, \qquad y=\eta.
\label{SWW 4}
\end{equation}
In this article we consider \emph{fully localised solitary waves}, that is nontrivial
travelling-wave solutions to \eqref{SWW 1}--\eqref{SWW 4} of the form
$\eta(x,z,t)=\eta(x-ct,z)$, $\varphi(x,y,z,t)=\varphi(x-ct,y,z)$
(so that the waves move with unchanging shape and constant speed $c$ from left to right)
with $\eta(x-ct,z) \rightarrow 0$ as $|(x-ct,z)| \rightarrow \infty$ (so that the waves
decay in every horizontal direction).

\begin{theorem}
Suppose that $c^2=2(1-\varepsilon^2)$. For each sufficiently small value of
$\varepsilon>0$ there exist two solitary-wave solutions
of \eqref{SWW 1}--\eqref{SWW 4} for which $\eta \in H^3(\R^2)$ is 
symmetric in $x$ and $z$ and given by
$$\eta(x,z)=\pm \varepsilon \zeta_0(\varepsilon x,\varepsilon z)\cos x + o(\varepsilon)$$
uniformly over $(x,z) \in \R^2$, where $\zeta_0$ is the unique
symmetric, positive (real) solution of the two-dimensional nonlinear Schr\"{o}dinger equation
\begin{equation}
-\tfrac{1}{2}\zeta_{xx}-\zeta_{zz} + \zeta -\tfrac{11}{16}|\zeta|^2 \zeta =0. \label{NLS SW}
\end{equation}
\end{theorem}

\begin{figure}[h]
\centering

\includegraphics[width=7.5cm]{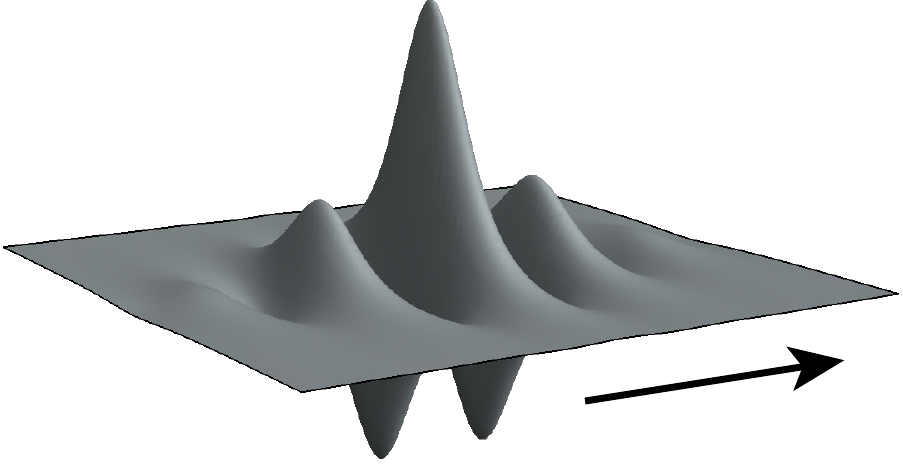}\hspace{0.5cm}\includegraphics[width=7cm]{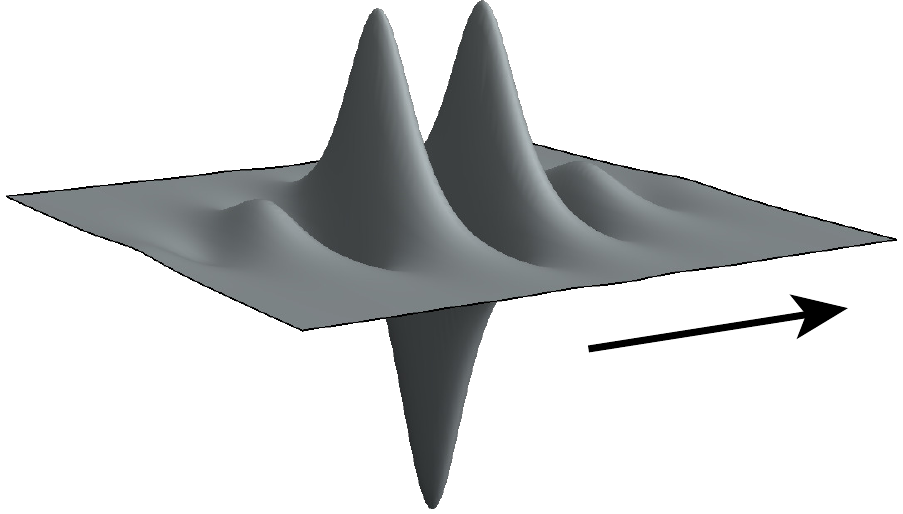}
{\it
\caption{Sketch of a symmetric fully localised solitary wave of elevation (left) and depression (right); the
arrow shows the direction of wave propagation. \label{flocs}}}
\end{figure}

This result confirms the prediction made on the basis of model equations (see below) and
numerical computations by Parau, Vanden-Broeck \& Cooker \cite{ParauVandenBroeckCooker05b}
(see Figure \ref{flocs} for sketches of typical free surfaces in their simulations). Qualitative
properties of (two- and three-dimensional) solitary waves on deep water have been discussed
by Wheeler \cite{Wheeler18}.

We proceed by formulating the water-wave problem \eqref{SWW 1}--\eqref{SWW 4}
in terms of the variables $\eta$ and $\Phi=\varphi|_{y=\eta}$ (see Zakharov \cite{Zakharov68} and
Craig \& Sulem \cite{CraigSulem93}). The Zakharov-Craig-Sulem formulation of the water-wave problem  is
\begin{align*}
& \eta_t - G(\eta)\Phi =0, \\
& \Phi_t +\eta
+ \frac{1}{2}\Phi_x^2+ \frac{1}{2}\Phi_z^2 -\frac{(G(\eta)\Phi+\eta_x\Phi_x+\eta_z\Phi_z)^2}{2(1+\eta_x^2+\eta_z^2)}
\!-\!\left[\frac{\eta_x}{\sqrt{1+\eta_x^2+\eta_z^2}}\right]_{\!x}
\!\!-\!\left[\frac{\eta_z}{\sqrt{1+\eta_x^2+\eta_z^2}}\right]_{\!z} \!\!\!=\! 0,
\end{align*}
where $G(\eta)\Phi=\varphi_y -\eta_x\varphi_x-\eta_z\varphi_z\big|_{y=\eta}$ and
$\varphi$ is the (unique) solution of the boundary-value problem
\begin{eqnarray*}
& & \parbox{6cm}{$\varphi_{xx}+\varphi_{yy}+\phi_{zz}=0,$}y<\eta, \\
& & \parbox{6cm}{$\varphi_y \rightarrow 0,$}y \rightarrow -\infty, \\
& & \parbox{6cm}{$\varphi=\Phi,$}y=\eta.
\end{eqnarray*}
Travelling waves are solutions of the form $\eta(x,z,t)=\eta(x-ct,z)$, $\Phi(x,z,t)=\Phi(x-ct,z)$; they satisfy
\begin{align}
& -c \eta_x - G(\eta)\Phi =0, \label{TW ZCS 1} \\
& -c\Phi_x +\eta
+ \frac{1}{2}\Phi_x^2+ \frac{1}{2}\Phi_z^2 \nonumber \\
& \hspace{1.5cm}\mbox{}-\frac{(G(\eta)\Phi+\eta_x\Phi_x+\eta_z\Phi_z)^2}{2(1+\eta_x^2+\eta_z^2)}
-\left[\frac{\eta_x}{\sqrt{1+\eta_x^2+\eta_z^2}}\right]_{x}
-\left[\frac{\eta_z}{\sqrt{1+\eta_x^2+\eta_z^2}}\right]_{z} = 0. \label{TW ZCS 2}
\end{align}

Equations \eqref{TW ZCS 1}, \eqref{TW ZCS 2} can be reduced to a single equation
for $\eta$. Using \eqref{TW ZCS 1}, one finds that $\Phi=-cG(\eta)^{-1}\eta_x$, and inserting this formula into
\eqref{TW ZCS 2} yields the equation
\begin{equation}
\KKp(\eta)-c^2\LLp(\eta)=0, \label{Basic equation - intro}
\end{equation}
where
\begin{align}
\KKp(\eta)&=\eta -\left[\frac{\eta_x}{\sqrt{1+\eta_x^2+\eta_z^2}}\right]_{\!x}
-\left[\frac{\eta_z}{\sqrt{1+\eta_x^2+\eta_z^2}}\right]_{\!z}, \label{Formula for KK} \\
\LLp(\eta) &=-\frac{1}{2}(K(\eta)\eta)^2-\frac{1}{2}(L(\eta)\eta)^2
+ \frac{(\eta_x-\eta_xK(\eta)\eta-\eta_zL(\eta)\eta)^2}{2(1+\eta_x^2+\eta_z^2)} + K(\eta)\eta
\label{Formula for LL}
\end{align}
and
$$K(\eta)\xi = -(G(\eta)^{-1} \xi_x)_x, \qquad L(\eta)\xi = -(G(\eta)^{-1} \xi_x)_z.$$
Note the equivalent definitions
\begin{equation}
K(\eta)\xi=-(\varphi|_{y=\eta})_x, \qquad L(\eta)\xi=-(\varphi|_{y=\eta})_z, \label{Formulae for K, L}
\end{equation}
where $\varphi$ is the solution of the boundary-value problem
\begin{eqnarray}
& & \parbox{6cm}{$\varphi_{xx}+\varphi_{yy}+\varphi_{zz}=0,$}y<\eta, \label{K 1}\\
& & \parbox{6cm}{$\varphi_y \rightarrow 0,$}y \rightarrow -\infty, \label{K 2}\\
& & \parbox{6cm}{$\varphi_y -\eta_x\varphi_x-\eta_z\varphi_z=\xi_x,$}y=\eta \label{K 3}
\end{eqnarray}
(which is unique up to an additive constant); the operators $K$ and $L$ are studied in Section \ref{Anal} below.
Although this fact is not used in the present paper, let us note that \eqref{Basic equation - intro} is in fact the Euler-Lagrange equation for the functional
$$
\JJ(\eta) := \KK(\eta) - c^2 \LL(\eta),
$$
where
$$
\KK(\eta)=\int_{\R^2} \left(\frac12 \eta^2 +\beta\sqrt{1+\eta_x^2+\eta_z^2}-\beta\right)\dx\dz,
\qquad
\LL(\eta)=\frac{1}{2}\int_{\R^2} \eta\, K(\eta) \eta \dx\dz;
$$
the functions $\KKp$ and $\LLp$ are the $L^2(\R^2)$-gradients of respectively $\KK$
and $\LL$ (see Buffoni \emph{et al.}\ 
\cite{BuffoniGrovesSunWahlen13,BuffoniGrovesWahlen18}).
Finally, observe that equation \eqref{Basic equation - intro} is invariant under the reflections
$\eta(x,z) \mapsto \eta(-x,z)$ and $\eta(x,z) \mapsto \eta(x,-z)$; a solution which is invariant
under these transformations is termed \emph{symmetric}.

\begin{figure}[h]
\centering

\includegraphics[scale=0.67]{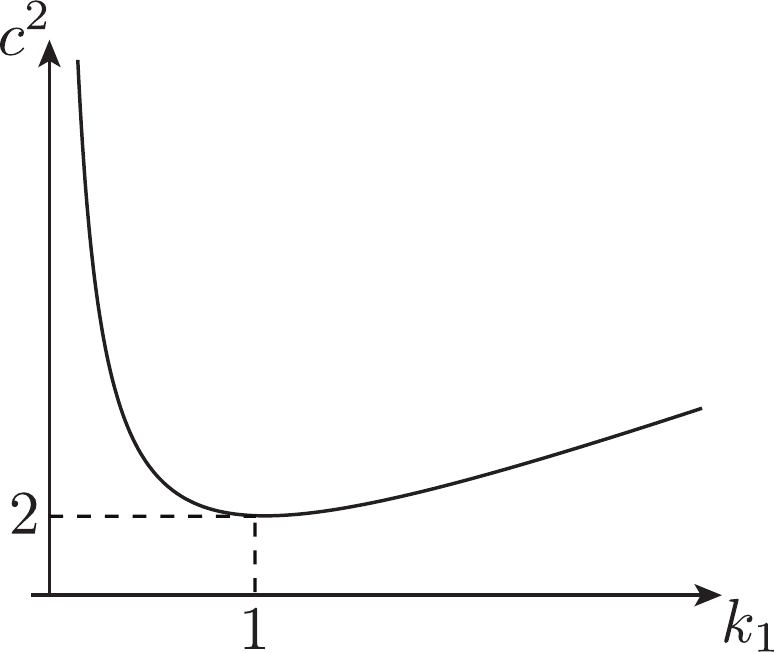}
{\it
\caption{Dispersion relation for a two-dimensional travelling wave train
with wave number $k_1 \geq 0$ and speed $c>0$.\label{Dispersion relation}}}
\end{figure}

It is instructive to review the formal derivation of the nonlinear Schr\"{o}dinger equation for travelling waves
(see Ablowitz \& Segur \cite[\S2.2]{AblowitzSegur79}),
beginning with sinusoidal wave trains.
The linearised version of \eqref{Basic equation - intro} admits a solution
of the form
$$\eta(x,z,t)=A \cos k_1(x-ct)$$
whenever $c>0$ and $k_1 \geq 0$ satisfy the linear dispersion relation
$$c^2 =k_1+\frac{1}{k_1}$$
(see Figure \ref{Dispersion relation}); note that the function $s \mapsto c(k_1)$,
$k_1 \geq 0$ has a unique global minimum $c_\mathrm{min}=\sqrt{2}$ at $k_1=1$.
Bifurcations of nonlinear solitary waves are expected whenever the
linear group and phase speeds are equal, so that $c^\prime(k_1)=0$ (see Dias \& Kharif \cite[\S 3]{DiasKharif99}).
We therefore expect the existence of small-amplitude solitary waves with speed near $\sqrt{2}$;
the waves bifurcate from a linear sinusoidal wave train with unit wavenumber.
Substituting $c^2=2(1-\varepsilon^2)$ and the \emph{Ansatz}
\begin{eqnarray*}
\lefteqn{\eta(x,z) = \tfrac{1}{2}\varepsilon \big(A_1(X,Z) \ee^{\ii  x}
+  \overline{A_1(X,Z)}\ee^{-\ii  x} \big)} \\
& & \hspace{0.75in}\mbox{}+\varepsilon^2 A_0(X,Z) +\tfrac{1}{2}\varepsilon^2 \big(A_2(X,Z) \ee^{2\ii  x}
+\overline{A_2(X,Z)}\ee^{-2\ii  x}\big) + \cdots,
\end{eqnarray*}
where $X=\varepsilon x$, $Z=\varepsilon z$,
into equation \eqref{Basic equation - intro}, one finds that $A_1$ satisfies the stationary nonlinear Schr\"{o}dinger equation
\eqref{NLS SW}.
This equation has a unique symmetric, positive (real) solution $\zeta_0 \in \SS(\R^2)$ which is characterised as the ground state of
the functional $\tilde{\JJ}: H^1(\R^2) \rightarrow {\mathbb C}$ with
$$\tilde{\JJ}(\zeta) = \int_{\R^2}\left(\frac{1}{4} |\zeta_x|^2
+\frac{1}{2} |\zeta_z|^2+\frac 1 2 |\zeta|^2
-\frac{11}{64}|\zeta|^4\right) \dx\dz$$
(see Sulem \& Sulem \cite[\S 4.2]{SulemSulem} and the references therein).

\begin{figure}[h]
\centering

\includegraphics[scale=0.7]{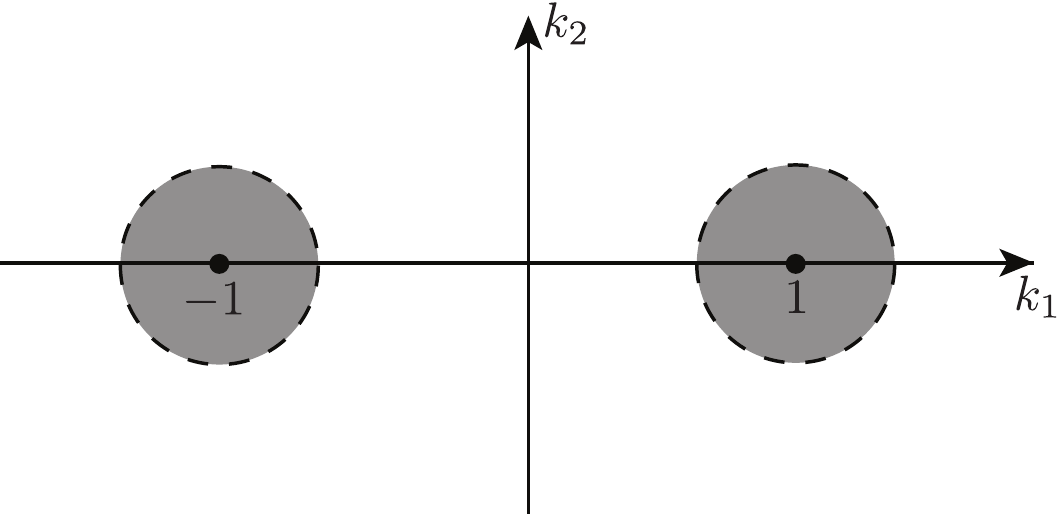}
{\it
\caption{The support of $\hat{\eta}_1$ is contained in the set $B=B_\delta(1,0) \cup B_\delta(-1,0)$. \label{Splitting}}}
\end{figure}

The above Ansatz suggests that the Fourier transform of a fully localised solitary wave 
is concentrated near the points $(1,0)$ and $(-1,0)$. We therefore decompose $\eta$ into the sum of functions
$\eta_1$ and $\eta_2$ whose Fourier transforms
$\hat{\eta}_1$ and $\hat{\eta}_2$ are supported in
the region $B=B_\delta(1,0) \cup B_\delta(-1,0)$ (with $\delta \in (0,\frac{1}{5})$) and its complement (see Figure \ref{Splitting}).
(The Fourier transform $\hat{u}=\FF[u]$ is defined by the formula
$$\hat{u}(k)=\frac{1}{2\pi}\int_{\R^2}u(x,z)\ee^{-\ii(k_1x + k_3z)}\dx\dz,\qquad k=(k_1,k_3),$$
and we use the notation $m(D)$ with $D=(-\mathrm{i}\partial_x,-\mathrm{i}\partial_z)$ for the Fourier multiplier-operator with symbol $m$, so that
$m(D)u = \FF^{-1}[m\hat{u}]$; in particular $\eta_1 = \chi(D)\eta$,
$\eta_2 = (1-\chi(D))\eta$, where $\chi$ is the characteristic function of the set $B$.)
Writing $c^2=2(1-\varepsilon^2)$ and
decomposing \eqref{Basic equation - intro} into
\begin{align*}
\chi(D)\big(\KKp(\eta_1+\eta_2)-2(1-\varepsilon^2)\LLp(\eta_1+\eta_2)\big)&=0, \\
(1-\chi(D))\big(\KKp(\eta_1+\eta_2)-2(1-\varepsilon^2)\LLp(\eta_1+\eta_2)\big)&=0,
\end{align*}
one finds that the second equation can be solved for $\eta_2$ as a function of $\eta_1$ for sufficiently
small values of $\varepsilon$; substituting 
$\eta_2=\eta_2(\eta_1)$ into the first yields the reduced equation
$$\chi(D)\big(\KKp(\eta_1+\eta_2(\eta_1))-2(1-\varepsilon^2)\LLp(\eta_1+\eta_2(\eta_1))\big)=0$$
for $\eta_1$. Finally, the scaling
\begin{equation}
\eta_1(x,z) = \tfrac{1}{2}\varepsilon \zeta(X,Z) \ee^{\ii x} + \tfrac{1}{2}\varepsilon \overline{\zeta(X,Z)}\ee^{-\ii x}
\label{NLS scaling as used here}
\end{equation}
transforms the reduced equation into a perturbation of the equation
\begin{equation}
\varepsilon^{-2}g(e+\varepsilon D)\zeta + 2f(e+\varepsilon D)\zeta - \tfrac{11}{8}|\zeta|^2\zeta=0, \label{Red eq - intro}
\end{equation}
where $e=(1,0)$ and
$$g(s)=1+|s|^2 -2f(s), \qquad f(s)=\frac{s_1^2}{|s|}$$
(see Sections \ref{sec:reduction} and \ref{sec:redeq}; the reduced equation is stated precisely in equation \eqref{PFDNLS}).
Equation \eqref{Red eq - intro} is termed a \emph{full-dispersion}
version of the stationary nonlinear Schr\"{o}dinger equation \eqref{NLS SW} since it retains the linear part of the original equation
\eqref{Basic equation - intro}; 
noting that
$$\varepsilon^{-2}g(e+\varepsilon k) + 2f(e+\varepsilon k) = 2+k_1^2+2k_3^2 + O(\varepsilon),$$
we obtain the fully reduced model equation
in the formal limit $\varepsilon=0$
(see Obrecht \& Saut \cite{ObrechtSaut15} for a discussion of related full-dispersion model equations for
three-dimensional water waves).

The existence theory is completed in Section \ref{sec:existence}, where we
exploit the fact that the reduction procedure preserves the invariance
of equation \eqref{Basic equation - intro} under
$\eta(x,z) \mapsto \eta(-x,z)$ and $\eta(x,z) \mapsto \eta(x,-z)$, so that equation \eqref{Red eq - intro} is
invariant under the reflections $\zeta(x,z) \mapsto \overline{\zeta(-x,z)}$ and $\zeta(x,z) \mapsto \zeta(x,-z)$.
We demonstrate that the reduced equation for $\zeta$ has two symmetric solutions $\zeta_\varepsilon^\pm$
which satisfy $\zeta_\varepsilon^\pm \rightarrow \pm \zeta_0$ in $H^1(\R^2)$ as $\varepsilon
\rightarrow 0$. The key step is a nondegeneracy result for the solution $\zeta_0$ of \eqref{NLS SW}
(see Weinstein \cite{Weinstein85}, Kwong \cite{Kwong89} and
Chang \emph{et al.}\ \cite{ChangGustafsonNakanishiTasi07}) in a symmetric setting
which allows one to apply a suitable version of the implicit-function theorem. 
A similar method was recently used by Stefanov \& Wright \cite{StefanovWright20} to establish the existence of
solitary-wave solutions to the Whitham equation (a full-dispersion Korteweg-de Vries equation).

The scaling \eqref{NLS scaling as used here} implies that our waves have small amplitude but finite energy.
When splitting our basic function space
$\XX=H^3(\R^2)$ into two parts $\XX_1=\chi(D)\XX$, $\XX_2=(1-\chi(D))\XX$ for $\eta_1$ and $\eta_2$,
we respect this scaling by equipping $\XX_1$ with the scaled norm $\nn\cdot\nn$ defined by
\begin{equation}
\nn \eta_1 \nn^2:= \int_{\R^2} (1+\varepsilon^{-2}((|k_1|-1)^2+k_3^2))|\hat \eta_1(k)|^2\dk_1 \dk_3 =\|\zeta\|_{H^1(\R^2)}^2 \label{Defn of triple norm}
\end{equation}
and taking $\zeta$ in a ball $B_R(0) \subseteq H^1(\R^2)$ which is large enough to contain $\zeta_0$; solving
the equation for $\eta_2$ yields the estimate
$$\|\eta_2(\eta_1)\|_{H^3(\R^2)} \lesssim \varepsilon^\theta \nn \eta_1 \nn^2,$$
where $\theta$ is a fixed number in the interval $(0,1)$. Equation \eqref{Defn of triple norm} shows that our waves have finite $H^3(\R^2)$-norm, while the estimates
$$\|\eta_1\|_\infty \lesssim \varepsilon^\theta \nn \eta_1 \nn, \qquad \|\eta_2\|_\infty \lesssim \|\eta_2\|_{H^3(\R^2)},$$
shows that they have small amplitude.

Our result complements recent existence theories for fully localised gravity-capillary solitary waves on
water of finite depth (Groves \& Sun \cite{GrovesSun08} and Buffoni \emph{et al.}\ 
\cite{BuffoniGrovesSunWahlen13,BuffoniGrovesWahlen18}), which also confirm predictions made
by model equations, namely the KP-I equation for `strong' and Davey-Stewartson equation
for `weak' surface tension (see Ablowitz \& Segur \cite{AblowitzSegur79}). In particular,
Buffoni \emph{et al.}\ \cite{BuffoniGrovesWahlen18} present a variational counterpart of the
theory in the present paper by reducing a classical variational principle for fully localised solitary waves
to a locally equivalent variational principle featuring a perturbation of the functional associated with the
Davey-Stewartson equation. A nontrivial critical point of the reduced functional is found by showing
that an appropriate direct method for the Davey-Stewartson functional (minimisation over its natural constraint set)
is robust under perturbation. This variational method is also applicable here, allowing one to reduce the
functional $\JJ$ to a perturbation of the functional $\tilde{\JJ}$. The present method however has the
advantages of being more explicit and yielding two distinct families of fully localised solitary waves.\pagebreak

\section{Analyticity} \label{Anal}

In this section we show that the operators $K$, $L$ given by \eqref{Formulae for K, L}
and hence $\KKp$ and $\LLp$ given by \eqref{Formula for KK}, \eqref{Formula for LL} are analytic
at the origin in suitable function spaces (see Corollaries \ref{K, L are anal} and \ref{KK, LL are anal} below).

The boundary-value problem \eqref{K 1}--\eqref{K 3} is handled using the change of variable
$$
y^\prime=y - \eta(x,z),\qquad u(x,y^\prime,z)=\varphi(x,y,z),
$$
which maps $\Sigma_\eta=\{(x,y,z)\colon x,z \in \R, -\infty < y < \eta(x,z)\}$ to the lower half-space
$\Sigma=\R \times (-\infty,0) \times \R$.
Dropping the primes, one finds that \eqref{K 1}--\eqref{K 3} are transformed into
\begin{align}
&u_{xx}+u_{yy}+u_{zz}=\partial_x F_1(\eta,u)+\partial_y F_2(\eta,u)+\partial_z F_3(\eta,u), && \quad y<0, \label{BC for u 1} \\
&u_y \rightarrow 0,&& \quad y \rightarrow - \infty, \label{BC for u 2} \\
&u_y=F_2(\eta,u)+\xi_x, && \quad y=0, \label{BC for u 3}
\end{align}
where
\[
F_1(\eta,u)=\eta_x u_y, \qquad
F_2(\eta,u)=\eta_x u_x + \eta_z u_z -(\eta_x^2+\eta_z^2) u_y, \qquad
F_3(\eta,u)=\eta_z u_y
\]
and $K(\eta)\xi = -u_x|_{y=0}$, $L(\eta)\xi=-u_z|_{y=0}$. We study this boundary-value problem in the space
$$\ZZ=\{\eta \in \SS^\prime(\R^2): \|\eta\|_\ZZ:=\|\hat{\eta}_1\|_{L^1(\R^2)}+\|\eta_2\|_3 < \infty\}$$
for $\eta$ and $H^3_\star(\Sigma)$ for $u$, in which
$H_\star^{n+1}(\Sigma)$,  $n \in {\mathbb N}$, is the completion of
$$
\SS(\Sigma,{\mathbb R})
=\{u \in C^\infty(\overline{\Sigma}):
|(x,z)|^m|\partial_x^{\alpha_1}\partial_y^{\alpha_2}\partial_z^{\alpha_3}u|\mbox{ is bounded for all }m,\alpha_1,\alpha_2, \alpha_3 \in {\mathbb N}_0 \}
$$
with respect to the norm
$$\|u\|_{n+1,\star}^2 := \|u_x\|_n^2 + \|u_y\|_n^2 + \|u_z\|_n^2$$
and $\|\cdot\|_s$ denotes the usual norm for the standard Sobolev space $H^s(\R^2)$ or $H^s(\Sigma)$.

\begin{lemma} \label{u is anal}
For each $\xi \in H^{5/2}(\R^2)$ and sufficiently small $\eta \in \ZZ$
the boundary-value problem \eqref{BC for u 1}--\eqref{BC for u 3} admits a unique solution
$u \in H^3_\star(\Sigma)$.
Furthermore, the mapping $\eta \mapsto (\xi \mapsto u)$ defines a function
$\ZZ\to \LL(H^{5/2}(\R^2), H^3_\star(\Sigma))$ which is analytic at the origin.
\end{lemma}
\proof
First note that for each $F_1$, $F_2$, $F_3 \in H^2(\Sigma)$ and $\xi \in H^{5/2}(\R^2)$ the boundary-value problem
\begin{align*}
&u_{xx}+u_{yy}+u_{zz}=\partial_x F_1+ \partial_y F_2 +\partial_z F_3, && \quad y<0, \\
&u_y \rightarrow 0,&& \quad y \rightarrow - \infty, \\
&u_y=F_2(\eta,u)+\xi_x, && \quad y=0,
\end{align*}
admits a unique solution $u=S(F_1,F_2,F_3,\xi)$ in $H^3_\star(\Sigma)$ whose gradient is obtained from the explicit formula
\begin{eqnarray*}
\lefteqn{S(F_1,F_2,F_3,\xi)=\FF^{-1}\left[\int_{-\infty}^0 \!\!\left(-\frac{\ii k_1}{2|k|}\hat{F}_1 - \frac{\ii k_3}{2|k|} \hat{F}_3 + \frac{1}{2}\sgn(y-\tilde{y}) \hat{F}_2\right)
\ee^{-|k||y-\tilde{y}|} \dtildey\right.}\hspace{1cm}\\
&& \hspace{1.25in}+ \left.\!\!\int_{-\infty}^0 \!\!\left(-\frac{\ii k_1}{2|k|}\hat{F}_1 - \frac{\ii k_3}{2|k|} \hat{F}_3+\frac{1}{2}\hat{F}_2\right)\ee^{|k|(y+\tilde{y})}\dtildey
+\frac{\ii k_1}{|k|}\hat{\xi}\ee^{|k|y}\right]
\end{eqnarray*}
(with a slight abuse of notation), so that
$$\|S(F_1,F_2,F_3,\xi)\|_{3,\star} \lesssim \|F_1\|_2 + \|F_2\|_2 + \|F_3\|_2 + \|\xi\|_{5/2}.$$

Define
$$T:H_\star^3(\Sigma) \times \ZZ \times H^{5/2}(\R^2)
\rightarrow H_\star^3(\Sigma)$$
by
$$T(u,\eta,\xi)=u-S(F_1(\eta,u),F_2(\eta,u),F_3(\eta,u),\xi)$$
and note that the solutions of \eqref{BC for u 1}--\eqref{BC for u 3} are precisely the zeros of $T(\cdot,\eta,\xi)$.
Using the estimates
\begin{align*}
\|\eta_x^j w\|_2 & \leq  \|\eta_{1x}^j w\|_2 + \|\eta_{2x}^j w\|_2\\
& \lesssim \|\eta_1\|_{3,\infty}^j \|w\|_2 + \|\eta_2\|_3^j \|w\|_2 \\
& \lesssim (\|\hat{\eta}_1\|_{L^1(\R^2)} + \|\eta_2\|_3)^j \|w\|_2 \\
&= \|\eta\|_\ZZ^j \|w\|_2, \qquad j=1,2,
\end{align*}
and similarly
$$\|\eta_z^j w\|_2 \lesssim \|\eta\|_\ZZ^j \|w\|_2, \qquad j=1,2$$
(where we have used the fact that $\hat{\eta}_1$ has compact support), we find
that the mappings $H_\star^3(\Sigma) \times \ZZ
\rightarrow H^2(\Sigma)$ given by
$(\eta,u) \mapsto F_j(\eta,u)$ , $j=1,2,3$,
are analytic at the origin; it follows that $T$ is also analytic at the origin.
Furthermore $T(0,0,0)=0$ and 
$\mathrm{d}_1T[0,0,0]=I$ is an isomorphism.
By the analytic implicit-function theorem there exist open neighbourhoods $V_1$ and $V_2$ of the origin in $\ZZ$
and $H^{5/2}(\R^2)$ and an analytic function $v: V_1 \times V_2 \rightarrow H_\star^3(\Sigma)$
such that
$$T(v(\eta,\xi),\eta,\xi)=0.$$
Since $v$ is linear in $\xi$ one can take $V_2$ to be the entire space $H^{5/2}(\R^2)$.\qed

\begin{corollary} \label{K, L are anal}
The mappings $K(\cdot), L(\cdot) \colon \ZZ \to \LL(H^{5/2}(\R^2), H^{3/2}(\R^2))$
are analytic at the origin.
\end{corollary}

In view of Corollary \ref{K, L are anal} we choose $M$ sufficiently small and study the equation
\begin{equation}
\KKp(\eta)-2(1-\varepsilon^2)\LLp(\eta)=0 \label{Basic equation}
\end{equation}
in the set
$$U=\{\eta \in H^3(\R^2): \|\eta\|_\ZZ < M\},$$
noting that $H^3(\R^2)$ is continuously embedded in $\ZZ$ and that $U$ is an open neighbourhood of the origin in $H^3(\R^2)$;
we proceed accordingly by decomposing $\XX=H^3(\R^2)$ into the direct sum of
$\XX_1=\chi(D)H^3(\R^2)$ and $\XX_2=(1-\chi(D))H^3(\R^2)$.

\begin{corollary} \label{KK, LL are anal}
The formulae \eqref{Formula for KK}, \eqref{Formula for LL}
define functions $U \rightarrow H^1(\R^2)$
which are analytic at the origin and satisfy $\KKp(0)=\LLp(0)=0$.
\end{corollary}
\proof The result for $\KKp$ follows from \eqref{Formula for KK}, Corollary \ref{K, L are anal} and the fact that
$H^{3/2}(\R^2)$ is an algebra. The result for $\LLp$ follows from \eqref{Formula for LL} and the observation that
$\eta \mapsto \eta_x (1+\eta_x^2+\eta_z^2)^{-1/2}$ and $\eta \mapsto \eta_z (1+\eta_x^2+\eta_z^2)^{-1/2}$
define functions $U \rightarrow H^2(\R^2)$ which are analytic at the origin since $H^2(\R^2)$ is an
algebra.\qed\medskip

In keeping with Lemma \ref{u is anal} and Corollary \ref{K, L are anal} we write
$$u(\eta,\xi) = \sum_{j=0}^\infty u^j(\eta,\xi),$$
where $u^j$ is homogeneous of degree $j$ in $\eta$ and linear in $\xi$, and
$$K(\eta)=\sum_{j=0}^\infty K_j(\eta), \quad L(\eta)=\sum_{j=0}^\infty L_j(\eta),
\quad
\KKp(\eta):=\sum_{j=1}^\infty \KKp_j(\eta), \quad \LLp(\eta):=\sum_{j=1}^\infty \LLp_j(\eta).$$
where $K_j(\eta)$, $L_j(\eta)$ and $\KKp_j(\eta)$, $\LLp_j(\eta)$ are homogeneous of degree $j$ in $\eta$.
A straightforward calculation shows that
$$u^0(\xi) = \FF^{-1}\left[\frac{\i k_1}{|k|}\ee^{|k|y}\hat{\xi}\right]$$
and hence that $K_0$ and $L_0$ are Fourier-multiplier operators, namely
$$K_0\xi = \FF^{-1}\left[\frac{k_1^2}{|k|} \hat{\xi}\right], \qquad L_0\xi = \FF^{-1}\left[\frac{k_1k_3}{|k|} \hat{\xi}\right]$$
(we have omitted the argument $\eta$ on the left-hand sides of these equations).

The following lemma gives
expressions for the first few terms in the Maclaurin expansions of $\KKp(\eta)$ and $\LLp(\eta)$;
it is proved by expanding \eqref{Formula for KK}, \eqref{Formula for LL} and
examining the boundary-value problems for $u^1(\eta,\eta)$ and $u^2(\eta,\eta)$ to derive
the formulae
\begin{align}
K_1(\eta)\eta &= -(\eta \eta_x)_x-K_0(\eta K_0 \eta)-L_0(\eta L_0 \eta), \label{Formula for K1} \\
L_1(\eta)\eta &= -(\eta \eta_x)_z-L_0(\eta K_0 \eta)-M_0(\eta L_0 \eta) \label{Formula for L1}
\end{align}
with a similar formula for $K_2(\eta_1)\eta_1$ (see Buffoni \emph{et al.}\ \cite[pp.\ 1032--1033]{BuffoniGrovesSunWahlen13} for details in a similar setting;
the restriction to $\eta_1$ is necessary to allow the use of higher-order derivatives in these expressions).

\begin{lemma} $ $ \label{KK, LL expansions}
\begin{list}{(\roman{count})}{\usecounter{count}}
\item
The identities 
\begin{align}
\KKp_1(\eta)&=\eta- \eta_{xx} - \eta_{zz}, \nonumber \\
\KKp_2(\eta)&=0, \nonumber \\
\KKp_3(\eta)&=\frac{1}{2} ((\eta_x^2+\eta_z^2)\eta_x)_x + \frac{1}{2}(\eta_x^2+\eta_z^2)\eta_z)_z \label{Formula for KKp3}
\end{align}
hold for each $\eta \in H^3(\R^2)$.
\item
The identities
\begin{align*}
\LLp_1(\eta) &= K_0\eta, \\
\LLp_2(\eta) &= \tfrac{1}{2}\left(\eta_x^2 - (K_0 \eta)^2 -(L_0 \eta)^2
-2(\eta_x\eta)_x-2K_0(\eta K_0\eta)-2L_0(\eta L_0\eta) \right)
\end{align*}
hold for each $\eta \in H^3(\R^2)$.
\item
The identity
\begin{align*}
\LLp_3(\eta_1)&=K_0\eta_1\, K_0(\eta_1 K_0\eta_1) +K_0\eta_1 \,L_0 (\eta_1 L_0 \eta_1)
+L_0\eta_1 \,L_0(\eta_1 K_0\eta_1)+L_0\eta_1 \, M_0(\eta_1 L_0 \eta_1)\\
& \qquad\mbox{}+K_0(\eta_1 K_0(\eta_1 K_0\eta_1)) +K_0(\eta_1 L_0 (\eta_1 L_0 \eta_1))
+L_0(\eta_1 L_0(\eta_1 K_0\eta_1))\\
&\qquad\mbox{}+L_0(\eta_1 M_0(\eta_1 L_0 \eta_1))+\eta_1 (K_0 \eta_1)\eta_{1xx}+\tfrac{1}{2}K_0(\eta_1^2\eta_{1xx}) + \tfrac{1}{2}(\eta_1^2 K_0 \eta_1)_{xx} \\
& \qquad \mbox{}+\eta_1 (L_0 \eta_1)\eta_{1xz}+\tfrac{1}{2}L_0(\eta_1^2\eta_{1xz}) + \tfrac{1}{2}(\eta_1^2 L_0 \eta_1)_{xz},
\end{align*}
where
$$M_0\xi = \FF^{-1}\left[\frac{k_3^2}{|k|} \hat{\xi}\right],$$
holds for each $\eta_1 \in \XX_1$ and more generally for any function $\eta_1$ whose Fourier transform has compact support.
\end{list}
\end{lemma}

Finally, we present some useful estimates for their cubic and higher-order parts of $\KKp(\eta)$ and $\LLp(\eta)$. The results for $\LL^\prime(\eta)$ are
established by substituting 
$$K(\eta)=\sum_{j=0}^2 K_j(\eta)+K_\mathrm{c}(\eta), \qquad L(\eta)=\sum_{j=0}^2 L_j(\eta)+L_\mathrm{c}(\eta)$$
into \eqref{Formula for LL} and estimating the resulting formulae for $\LL_\mathrm{c}(\eta)$ and
$\LL_\mathrm{r}(\eta)$ using the rules
$$
\|K_j(\eta)\eta\|_{3/2} \lesssim \|\eta\|_\ZZ^j\|\eta\|_{5/2}, \qquad
\|K_\mathrm{c}(\eta)\eta\|_{3/2} \lesssim \|\eta\|_\ZZ^3\|\eta\|_{5/2}$$
(with corresponding estimates for $L_j(\eta)(\eta)$, $L_\mathrm{c}(\eta)(\eta)$ and derivatives). Since this method yields only
$$\|(K_1(\eta)\eta)^2\|_1,\ \|(L_1(\eta)\eta)^2\|_1 \lesssim \|\eta\|_\ZZ^2 \|\eta\|_3^2$$
we do not include the fourth-order terms $-\tfrac{1}{2}(K_1(\eta)\eta)^2$, $-\tfrac{1}{2}(L_1(\eta)\eta)^2$
in $\LLp_\mathrm{r}(\eta)$ and treat them separately later (see in particular Proposition \ref{Extra terms}).

\begin{lemma} $ $ \label{KK, LL hot estimates}
\begin{itemize}
\item[(i)]
The quantities
$$\KKp_\mathrm{c}(\eta):=\sum_{j=3}^\infty \KKp_j(\eta), \qquad \LLp_\mathrm{c}(\eta):=\sum_{j=3}^\infty \LLp_j(\eta)$$
satisfy the estimates
$$
\|\KKp_\mathrm{c}(\eta)\|_1 \lesssim \|\eta\|_\ZZ^2 \|\eta\|_3, \qquad
\|\mathrm{d}\KKp_\mathrm{c}[\eta](v)\|_1 \lesssim \|\eta\|_\ZZ^2 \|v\|_3 + \|\eta\|_\ZZ \|\eta\|_3 \|v\|_\ZZ
$$
$$
\|\LLp_\mathrm{c}(\eta)\|_1 \lesssim \|\eta\|_\ZZ^2 \|\eta\|_3, \qquad
\|\mathrm{d}\LLp_\mathrm{c}[\eta](v)\|_1 \lesssim \|\eta\|_\ZZ^2 \|v\|_3 + \|\eta\|_\ZZ \|\eta\|_3 \|v\|_\ZZ
$$
for each $\eta \in U$ and $v \in H^2(\R)$.
\item[(ii)]
The quantities
$$\KKp_\mathrm{r}(\eta):=\sum_{j=4}^\infty \KKp_j(\eta), \qquad \LLp_\mathrm{r}(\eta):=\sum_{j=4}^\infty \LLp_j(\eta)+\tfrac{1}{2}(K_1(\eta)\eta)^2+\tfrac{1}{2}(L_1(\eta)\eta)^2$$
satisfy the estimates
$$
\|\KKp_\mathrm{r}(\eta)\|_1 \lesssim \|\eta\|_{\ZZ}^4 \|\eta\|_3, \qquad
\|\mathrm{d}\KKp_\mathrm{r}[\eta](v)\|_1 \lesssim \|\eta\|_{\ZZ}^4 \|v\|_3 + \|\eta\|_{\ZZ}^3  \|\eta\|_3\|v\|_\ZZ
$$
$$
\|\LLp_\mathrm{r}(\eta)\|_1 \lesssim \|\eta\|_{\ZZ}^3 \|\eta\|_3, \qquad
\|\mathrm{d}\LLp_\mathrm{r}[\eta](v)\|_1 \lesssim \|\eta\|_{\ZZ}^3 \|v\|_3 + \|\eta\|_{\ZZ}^2  \|\eta\|_3\|v\|_\ZZ
$$
for each $\eta \in U$ and $v \in H^2(\R)$.
\end{itemize}
\end{lemma}

\section{Reduction} \label{sec:reduction}

Observe that
$\eta \in U$ satisfies \eqref{Basic equation} if and only if
\begin{align}
\eta_1-\eta_{1xx}-\eta_{1zz}-2 K_0\eta_1
+2\varepsilon^2K_0\eta_1+\chi(D)\NN(\eta_1+\eta_2)&=0,
\label{Pre X_1 component}  \\
\eta_2-\eta_{2xx}-\eta_{2zz}-2 K_0\eta_2 +2\varepsilon^2 K_0\eta_2+(1-\chi(D))\NN(\eta_1+\eta_2)&=0,
\label{X_2 component}
\end{align}
in which
$$\NN(\eta)=\KKp_\mathrm{c}(\eta)-2(1-\varepsilon^2)\big(\LLp_2(\eta)+\LLp_\mathrm{c}(\eta)\big).$$
The nonlinear term in \eqref{Pre X_1 component} is at leading order cubic in $\eta_1$
because $\chi(D)\LLp_2(\eta_1)$ vanishes; we therefore write it as
\begin{equation}
\eta_1-\eta_{1xx}-\eta_{1zz}-2 K_0\eta_1
+2\varepsilon^2K_0\eta_1+\chi(D)\left(\NN(\eta_1+\eta_2)+2(1-\varepsilon^2)\LLp_2(\eta_1)\right)=0
\label{X_1 component}
\end{equation}
and make the corresponding adjustment to \eqref{X_2 component}, that is `replacing' its nonlinearity with
$$(1-\chi(D))\left(\NN(\eta_1+\eta_2)+2(1-\varepsilon^2)\LLp_2(\eta_1)\right),$$
by writing
$$\eta_2 = F(\eta_1)+\eta_3, \qquad F(\eta_1):=
2(1-\varepsilon^2)\FF^{-1}\left[\frac{1-\chi(k)}{g(k)}\FF[\LLp_2(\eta_1)]\right]$$
(with the requirement that $\eta_1+F(\eta_1)+\eta_3 \in U$).
Equation \eqref{X_2 component} may thus be cast in the form
\begin{equation}
\eta_3 = -\FF^{-1}\Bigg[\frac{1-\chi(k)}{g(k)}\FF\Big[2(1-\varepsilon^2)\LLp_2(\eta_1)
+\NN(\eta_1+F(\eta_1)+\eta_3)+2\varepsilon^2 K_0(F(\eta_1)+\eta_3)\Big]\Bigg], \label{eta3 eqn}
\end{equation}
where
$$g(k)=1+|k|^2-2\frac{k_1^2}{|k|} \geq 0$$
with equality if and only if $k=\pm(1,0)$.

\begin{proposition}
\label{prop:bounded mapping}
The mapping $$f\mapsto \FF^{-1}\left[\frac{1-\chi(k)}{g(k)}\hat f\right]$$ defines a bounded linear operator $H^1(\R^2)\to H^3(\R^2)$.
\end{proposition}

We proceed by solving \eqref{eta3 eqn} for $\eta_3$ as a function of $\eta_1$ using the following following fixed-point theorem, which is a straightforward extension of a standard result in nonlinear analysis.

\begin{theorem}
\label{thm:fixed-point}
Let $\XX_1$, $\XX_2$ be Banach spaces, $X_1$, $X_2$ be closed, convex sets in, respectively, $\XX_1$, $\XX_2$ containing the origin and $\GG\colon X_1\times X_2 \to \XX_2$ be a smooth function. Suppose that there exists a 
continuous function $r\colon X_1\to [0,\infty)$ such that
$$
\|\GG(x_1,0)\|\le \tfrac{1}{2}r, \quad \|\mathrm{d}_2 \GG[x_1,x_2]\|\le \tfrac{1}{3}
$$ 
for each $x_2\in \overline B_r(0)\subseteq X_2$ and each $x_1\in X_1$.

Under these hypotheses there exists for each $x_1\in X_1$ a unique solution $x_2=x_2(x_1)$ of the fixed-point equation
$x_2=\GG(x_1,x_2)$
satisfying $x_2(x_1)\in \overline B_r(0)$. Moreover $x_2(x_1)$ is a smooth function of $x_1\in X_1$ and in particular satisfies the estimate
$$
\|\mathrm{d} x_2[x_1]\|\le 2\|\mathrm{d}_1 \GG[x_1, x_2(x_1)]\|.
$$
\end{theorem}

We apply Theorem \ref{thm:fixed-point} to equation \eqref{eta3 eqn} with
$\XX_1=\chi(D)H^3(\R^2)$, $\XX_2=(1-\chi(D))H^3(\R^2)$, equipping
$\XX_1$ with the scaled norm
$$
\nn \eta \nn:=\left( \int_{\R^2} (1+\varepsilon^{-2}((|k_1|-1)^2+k_3^2))|\hat \eta(k)|^2\dk_1 \dk_3\right)^{\!1/2}
$$
and $\XX_2$ with the usual norm for $H^3(\R^2)$, and taking
$$
X_1=\{\eta_1\in \XX_1 \colon \nn \eta_1\nn \le R_1\}, \qquad
X_3=\{\eta_3\in \XX_2 \colon \| \eta_3\|_3 \le R_3\};
$$
the function $\GG$ is given by the right-hand side of \eqref{eta3 eqn}. (Here we write
$X_3$ rather than $X_2$ for notational clarity.)
The calculation
\begin{align*}
\hspace{1.2cm}\int_{\R^2} |\hat \eta_1(k)|\dk_1\dk_3
&= \int_{\R^2}\frac{(1+\varepsilon^{-2}((|k_1|-1)^2+k_3^2))^{1/2}}{(1+\varepsilon^{-2}((|k_1|-1)^2+k_3^2))^{1/2}} |\hat \eta_1(k)| \dk_1 \dk_3\\
&\le 2\nn \eta\nn \left( \int_{B_\delta(1,0)} \frac{1}{1+\varepsilon^{-2}((k_1-1)^2+k_3^2)} \dk_1 \dk_3 \right)^{1/2} \\
& = 2\sqrt{\pi}\varepsilon (\log(1+\delta^2\varepsilon^{-2}))^{1/2} \nn \eta\nn
\end{align*}
shows that
\begin{equation}
\|\hat{\eta}_1\|_{L^1(\R^2)} \lesssim \varepsilon^\theta \nn \eta_1 \nn, \qquad \eta_1 \in \XX_1,
\label{Triple norm large}
\end{equation}
for each fixed $\theta \in (0,1)$. We can therefore guarantee that $\|\hat{\eta}_1\|_{L^1(\R^2)} < M/2$ for all
$\eta_1 \in X_1$ for an arbitrarily large value
of $R_1$; the value of $R_3$
is then constrained by the requirement that $\|F(\eta_1) + \eta_3\|_3 < M/2$ for all $\eta_1 \in X_1$ and $\eta_3 \in X_3$,
so that $\eta_1+F(\eta_1)+\eta_3 \in U$
(Corollary \ref{cor:F estimates} below asserts that $\|F(\eta_1)\|_3 = O(\varepsilon^\theta)$ uniformly over $\eta_1 \in X_1$).

We proceed by systematically estimating each term
appearing in the equation for $\GG$, using the inequalities
$$\|\eta\|_\infty \lesssim \|\eta\|_\ZZ, \quad \|\eta\|_\ZZ \lesssim \varepsilon^\theta \nn \eta_1 \nn + \|\eta_3\|_3, \quad
\|\eta\|_3 \lesssim \nn \eta_1 \nn + \|\eta_3\|_3$$
and making extensive use of the fact that the support of $\hat{\eta}_1$ is contained in the
fixed bounded set $B$, so that for example
$$\|\eta_1\|_n \lesssim \|\eta_1\|_0, \qquad \|\eta_1\|_{n,\infty} \lesssim \varepsilon^\theta \nn \eta_1 \nn$$
for each $n \in {\mathbb N}_0$.

In order to estimate $F(\eta_1)$ we write $\LLp_2(\eta)=m(\{\eta\}^2)$, where
\begin{align}
m(u,v)&=\tfrac{1}{2}\left(u_x v_x - (K_0 u) (K_0 v) -(L_0 u) (L_0 v) \right) \nonumber \\
&\qquad\mbox{}+\tfrac{1}{2}\left(
-(u_xv +u v_x)_x-K_0(u K_0 v+v K_0 u)-L_0(u L_0v +vL_0u)
\right) \label{Formula for m}
\end{align}
(see Lemma \ref{KK, LL expansions}(ii)), and note that
$$\mathrm{d}\LLp_2[\eta](v)=2m(\eta,v).$$

\begin{proposition}
\label{prop:m estimate}
The estimate
$$\|m(u,v)\|_1 \lesssim \|u\|_\ZZ \|v\|_3$$
holds for each $u$, $v\in H^3(\R^2)$.
\end{proposition}
\begin{corollary}
\label{cor:F estimates}
The estimates
$$\|F(\eta_1)\|_3\lesssim \varepsilon^\theta \nn \eta_1\nn^2,\quad
\|\mathrm{d}F[\eta_1]\|_{\LL(\XX_1, \XX_2)}  \lesssim \varepsilon^\theta \nn \eta_1\nn$$
hold for each $\eta_1\in X_1$.
\end{corollary}

\begin{remark} \label{rem:K0Feta1 remark}
Noting that
$$K_0 F(\eta_1) = 2(1-\varepsilon^2)\FF^{-1}\left[\frac{1-\chi(k)}{g(k)} \frac{k_1^2}{|k|} \FF[\LLp_2(\eta_1)]\right]$$
and that $\FF[\LLp_2(\eta_1)]$ has compact support, one finds
that $K_0 F(\eta_1)$ satisfies the same estimates as $F(\eta_1)$.
\end{remark}

\begin{lemma}
\label{lem:AA estimates}
The quantity
$$
\NN_1(\eta_1,\eta_3) = \LLp_2(\eta_1+F(\eta_1)+\eta_3)-\LLp_2(\eta_1)
$$
satisfies the estimates
\begin{list}{(\roman{count})}{\usecounter{count}}
\item
$\|\NN_1(\eta_1,\eta_3)\|_1\lesssim \varepsilon^{2\theta} \nn \eta_1\nn^3+\varepsilon^\theta \nn \eta_1\nn^2\|\eta_3\|_3
+\varepsilon^\theta \nn \eta_1\nn\|\eta_3\|_3+\|\eta_3\|_3^2$,
\item
$\|\mathrm{d}_1\NN_1[\eta_1,\eta_3]\|_{\LL(\XX_1,H^1(\R^2))}\lesssim \varepsilon^{2\theta} \nn \eta_1\nn^2
+\varepsilon^\theta \nn \eta_1\nn\|\eta_3\|_3+\varepsilon^\theta \|\eta_3\|_3$,
\item
$\|\mathrm{d}_2\NN_1[\eta_1,\eta_3]\|_{\LL(\XX_2,H^1(\R^2))}\lesssim \varepsilon^\theta \nn \eta_1\nn+\|\eta_3\|_3$
\end{list}
for each $\eta_1\in X_1$ and $\eta_3\in X_3$.
\end{lemma}
\proof
We estimate
$$\NN_1(\eta_1,\eta_3) = 2m(\eta_1, F(\eta_1)+\eta_3)+m(F(\eta_1)+\eta_3,F(\eta_1)+\eta_3)$$
by combining Proposition \ref{prop:m estimate} with Corollary \ref{cor:F estimates} using the chain rule.\hfill$\Box$

\begin{lemma} \label{BB estimates}
The quantity
$$
\NN_2(\eta_1,\eta_3)=\KKp_\mathrm{c}(\eta_1+F(\eta_1)+\eta_3)-2(1-\varepsilon^2)\LLp_\mathrm{c}(\eta_1+F(\eta_1)+\eta_3),
$$
satisfies the estimates
\begin{list}{(\roman{count})}{\usecounter{count}}
\item
$\|\NN_2(\eta_1,\eta_3)\|_1\lesssim (\varepsilon^\theta \nn \eta_1\nn+\|\eta_3\|_3)^2(\nn \eta_1\nn+\|\eta_3\|_3)$,
\item
$\|\mathrm{d}_1\NN_2[\eta_1,\eta_3]\|_{\LL(\XX_1,H^1(\R^2))} \lesssim (\varepsilon^\theta \nn \eta_1\nn+\|\eta_3\|_3)^2$,
\item
$\|\mathrm{d}_2\NN_2[\eta_1,\eta_3]\|_{\LL(\XX_2,H^1(\R^2))}\lesssim (\varepsilon^{\theta}\nn \eta_1\nn+\|\eta_3\|_3)(\nn \eta_1\nn+\|\eta_3\|_3)$
\end{list}
for each $\eta_1\in X_1$ and $\eta_3\in X_3$.
\end{lemma}
\proof We compute the derivatives of $\NN_2$ using the chain rule and estimate these
expressions using the linearity of the derivative,
Lemma \ref{KK, LL hot estimates}(i) and Corollary \ref{cor:F estimates}.
\qed\medskip

Altogether we have established the following estimates for $\GG$ and its derivatives
(see Proposition \ref{prop:bounded mapping},
Remark \ref{rem:K0Feta1 remark} and Lemmata \ref{lem:AA estimates}, \ref{BB estimates}).

\begin{lemma} \label{lem:Complete GG estimates}
The function $\GG: X_1 \times X_3 \rightarrow \XX_2$ satisfies the estimates
\begin{list}{(\roman{count})}{\usecounter{count}}
\item
$\|\GG(\eta_1,\eta_3)\|_3\lesssim (\varepsilon^\theta \nn \eta_1\nn+\|\eta_3\|_3)^2(1+\nn \eta_1\nn+\|\eta_3\|_3)+\varepsilon^2\|\eta_3\|_3$,
\item
$\|\mathrm{d}_1\GG[\eta_1,\eta_3]\|_{\LL(\XX_1,\XX_2)} \lesssim (\varepsilon^\theta \nn \eta_1\nn+\|\eta_3\|_3)
(\varepsilon^\theta+\varepsilon^\theta \nn \eta_1\nn+\|\eta_3\|_3)$,
\item
$\|\mathrm{d}_2\GG[\eta_1,\eta_3]\|_{\LL(\XX_2)}\lesssim (\varepsilon^{\theta}\nn \eta_1\nn+\|\eta_3\|_3)(1+\nn \eta_1\nn+\|\eta_3\|_3)+\varepsilon^2$
\end{list}
for each $\eta_1\in X_1$ and $\eta_3\in X_3$.
\end{lemma}

\begin{theorem} \label{thm:estimate eta3}
Equation \eqref{eta3 eqn} has a unique solution $\eta_3 \in
X_3$ which depends smoothly upon $\eta_1 \in X_1$ and satisfies the estimates
$$
\|\eta_3(\eta_1)\|_3 \lesssim \varepsilon^{2\theta} \nn \eta_1\nn^2, \qquad
\|\mathrm{d}\eta_3[\eta_1]\|_{\LL(\XX_1,\XX_2)} \lesssim \varepsilon^{2\theta} \nn \eta_1\nn.
$$
\end{theorem}
\proof
Choosing $R_3$ and $\varepsilon$ sufficiently small and setting
$r(\eta_1)=\sigma \varepsilon^{2\theta} \nn \eta_1 \nn^2$ for a sufficiently large value
of $\sigma>0$, one finds that
\[\|\GG(\eta_1,0)\|_3 \lesssim \tfrac{1}{2}r(\eta_1), \qquad
\|\mathrm{d}_2 \GG[\eta_1,\eta_3]\|_{\LL(\XX_2)}  \lesssim \varepsilon^{\theta}
\]
for $\eta_1 \in X_1$ and $\eta_3 \in \overline{B}_{r(\eta_1)}(0) \subset X_3$
(Lemma \ref{lem:Complete GG estimates}(i), (iii)).
Theorem \ref{thm:fixed-point} asserts that equation \eqref{eta3 eqn} has a unique solution $\eta_3$
in $\overline{B}_{r(\eta_1)}(0) \subset X_3$ which depends smoothly upon $\eta_1 \in X_1$, and
the estimate for its derivative follows from Lemma \ref{lem:Complete GG estimates}(ii).\qed\medskip

Substituting $\eta_2=\eta_1+F(\eta_1)+\eta_3(\eta_1)$ into \eqref{X_1 component} yields the reduced equation
\begin{align}
\eta_1-\eta_{1xx}&-\eta_{1zz}-2 K_0\eta_1 \nonumber \\
& \mbox{}+2\varepsilon^2K_0\eta_1+\chi(D)\left(\NN(\eta_1+F(\eta_1)+\eta_3(\eta_2))+2(1-\varepsilon^2)\LLp_2(\eta_1)\right)=0
\label{Red eq v1}
\end{align}
for $\eta_1 \in X_1$. Observe that this equation is invariant under the reflections $\eta_1(x,z) \mapsto \eta_1(-x,z)$
and $\eta_1(x,z) \mapsto \eta_1(x,-z)$;
a familiar argument shows that they are inherited from the corresponding invariance of \eqref{X_1 component}, \eqref{eta3 eqn} under $\eta_1(x,z) \mapsto \eta_1(-x,z)$, $\eta_3(x,z) \mapsto \eta_3(-x,z)$
and $\eta_1(x,z) \mapsto \eta_1(x,-z)$, $\eta_3(x,z) \mapsto \eta_3(x,-z)$
when applying Theorem \ref{thm:fixed-point}.

\section{Derivation of the reduced equation} \label{sec:redeq}

In this section we compute the leading-order terms in the reduced equation \eqref{Red eq v1}.
To this end we write
$$\eta_1 = \eta_1^+ + \eta_1^-,$$
where $\eta_1^+=\chi^+(D)\eta_1$, $\eta_1^-=\chi^-(D)\eta_1$
and $\chi^+$, $\chi^-$ are the characteristic functions of respectively $B_\delta(1,0)$ and
$B_\delta(-1,0)$, so that $\eta_1^+$ satisfies the equation
\begin{align}
\eta_1^+-\eta_{1xx}^+&-\eta_{1zz}^+-2 K_0\eta_1^+ \nonumber \\
&\mbox{}
+2\varepsilon^2K_0\eta_1+\chi^+(D)\left(\NN(\eta_1+F(\eta_1)+\eta_3(\eta_2))+2(1-\varepsilon^2)\LLp_2(\eta_1)\right)=0
\label{Reduced equation v2}
\end{align}
(and $\eta_1^-$ satisfies its complex conjugate). It is also convenient to introduce some additional notation.

\begin{definition} \hspace{1cm}
\begin{itemize}
\item[(i)]
The symbol $\underline{O}(\varepsilon^\gamma  \nn \eta_1 \nn^r)$ denotes a smooth function
$N: X_1 \rightarrow H^1(\R^2)$ which satisfies the estimates 
$$
\|N(\eta_1)\|_1 \lesssim \varepsilon^\gamma  \nn \eta_1 \nn^r, \qquad
\|\mathrm{d}N[\eta_1]\|_{\LL(\XX_1, H^1(\R^2))}\lesssim \varepsilon^\gamma  \nn \eta_1 \nn^{r-1}
$$
for each $\eta_1 \in X_1$ (where $\gamma \geq 0$, $r \geq 1$). Furthermore
$$\underline{O}_0(\varepsilon^\gamma  \nn \eta_1 \nn^r):=\chi_0(D)\underline{O}(\varepsilon^\gamma  \nn \eta_1 \nn^r), \qquad
\underline{O}_+(\varepsilon^\gamma  \nn \eta_1 \nn^r) := \chi^+(D)\underline{O}(\varepsilon^\gamma  \nn \eta_1 \nn^r),$$
where $\chi_0$ and $\chi^+$ are the characteristic functions of the sets $B_\delta(0,0)$ and
$B_\delta(1,0)$.
\item[(ii)]
The symbol $\underline{O}^\varepsilon_n(\| u \|_1^r)$
denotes $\chi_0(\varepsilon D)N(u)$, where $N$ is a smooth function\linebreak
$B_R(0) \subseteq \chi_0(\varepsilon D)H^1(\R^2) \rightarrow H^n(\R^2)$
or $B_R(0) \subseteq H^1(\R^2) \rightarrow H^n(\R^2)$ which satisfies the estimates
$$
\|N(u)\|_n \lesssim \| u \|_1^r, \quad
\|\mathrm{d}N[u]\|_{\LL(H^1(\R^2),H^n(\R^2))}\lesssim \| u \|_1^{r-1}
$$
for each $u \in B_R(0)$  (with $r \geq 1$, $n \geq 0$).
\end{itemize}
\end{definition}

We begin with a result which shows how a Fourier-multiplier operator $m(D)$ may be approximated by $m(\omega,0)$ when acting upon
a function whose Fourier transform is supported near the point $(\omega,0)$. Its proof is given by Buffoni \emph{et al.}\ \cite[Lemma 11]{BuffoniGrovesWahlen18} (in a slightly different context).

\begin{lemma}  \label{lem:approximate identities}
The estimates
\begin{list}{(\roman{count})}{\usecounter{count}}
\item
$\partial_x \eta_1^\pm = \pm \mathrm{i}\eta_1^\pm + \underline{O}(\varepsilon \nn \eta_1\nn)$,
\item
$\partial_x^2 \eta_1^\pm =-\eta_1^\pm+  \underline{O}(\varepsilon \nn \eta_1\nn)$,
\item
$\partial_z \eta_1^\pm =  \underline{O}(\varepsilon \nn \eta_1\nn)$,
\item
$K_0 \eta_1^\pm = \eta_1^\pm + \underline{O}(\varepsilon\nn \eta_1\nn)$,
\item
$L_0 \eta_1^\pm = \underline{O}(\varepsilon \nn \eta_1\nn)$,
\item
$K_0((\eta_1^\pm)^2) = 2(\eta_1^\pm)^2 +  \underline{O}(\varepsilon^{1+\theta}\nn \eta_1\nn^2)$,
 \item
$L_0((\eta_1^\pm)^2) =  \underline{O}(\varepsilon^{1+\theta}\nn \eta_1\nn^2)$,
\item
$K_0 (\eta_1^+ \eta_1^-) =  \underline{O}(\varepsilon^{1+\theta}\nn \eta_1\nn^2)$,
\item
$L_0 (\eta_1^+ \eta_1^-) =  \underline{O}(\varepsilon^{1+\theta}\nn \eta_1\nn^2)$,
\item
$\FF^{-1}[g(k)^{-1}\FF[(\eta_1^\pm)^2]]=(\eta_1^\pm)^2+ \underline{O}(\varepsilon^{1+\theta}\nn \eta_1\nn^2)$,
\item
$K_0(\eta_1^-(\eta_1^+)^2) = \eta_1^-(\eta_1^+)^2 +\underline{O}(\varepsilon^{1+2\theta}\nn \eta_1\nn^3)$,
\end{list}
hold for each $\eta_1 \in X_1$.
\end{lemma}

We proceed by approximating each term
in the nonlinearity on the right-hand side of \eqref{Reduced equation v2} according
to the rules given in Lemma \ref{lem:approximate identities}.

\begin{proposition} \label{pm expansion of F}
The estimate
$$
F(\eta_1)=-2\left((\eta_1^+)^2+(\eta_1^-)^2\right)+F_\mathrm{r}(\eta_1),\qquad F_\mathrm{r}(\eta_1)=
\underline{O}(\varepsilon^{1+\theta}\nn \eta_1\nn^2)
$$
holds for each $\eta_1 \in X_1$.
\end{proposition}
\proof
Using the expansions given in Lemma \ref{lem:approximate identities},
we find that
$$\LLp_2(\eta_1)=m(\eta_1,\eta_1)=-\left((\eta_1^+)^2+(\eta_1^-)^2\right)+\underline{O}(\varepsilon^{1+\theta}\nn \eta_1\nn^2).$$
It follows that
$$
\FF^{-1}\left[\frac{1-\chi(k)}{g(k)}\FF[\LLp_2(\eta_1)]\right] 
= - \left((\eta_1^+)^2+(\eta_1^-)^2\right)+\underline{O}(\varepsilon^{1+\theta}\nn \eta_1\nn^2)
$$
because of Lemma \ref{lem:approximate identities}(x) and the fact that
$$
\FF^{-1}\left[\frac{1-\chi(k)}{g(k)}\FF[\underline{O}(\varepsilon^{1+\theta} \nn \eta_1\nn^2)]\right] = \underline{O}(\varepsilon^{1+\theta} \nn \eta_1\nn^2)
$$
(because  $(1-\chi(k))g(k)^{-1}$ is bounded). We conclude that
$$F(\eta_1) = 2(1-\varepsilon^2)\FF^{-1}\left[\frac{1-\chi(k)}{g(k)}\FF[\LLp_2(\eta_1)]\right] =
-2\left((\eta_1^+)^2+(\eta_1^-)^2\right)+\underline{O}(\varepsilon^{1+\theta}\nn \eta_1\nn^2).\eqno{\Box}$$

\begin{remark} \label{Fr estimate}
The remainder term $F_\mathrm{r}(\eta_1)$ in the formula for $F(\eta_1)$ given in Proposition \ref{pm expansion of F}
satisfies
$$
\|F_\mathrm{r}(\eta_1)\|_n \lesssim \varepsilon^{1+\theta}  \nn \eta_1 \nn^2, \qquad
\|\mathrm{d}F_\mathrm{r}[\eta_1]\|_{\LL(\XX_1, H^n(\R^2))}\lesssim \varepsilon^{1+\theta}  \nn \eta_1 \nn
$$
for all $n \in \N_0$ since its Fourier transform is supported in the region $B+B$.
\end{remark}

\begin{proposition} \label{Compute redeq 1}
The estimate
$$\chi^+(D)\NN_1(\eta_1,\eta_3(\eta_1))= 4\chi^+(D)\big(\eta_1^-(\eta_1^+)^2\big) + \underline{O}_+(\varepsilon^{3\theta}\nn \eta_1 \nn^3)$$
holds for each $\eta_1 \in X_1$.
\end{proposition}
\proof
Observe that
\begin{align*}
\chi^+(D)\NN_1(\eta_1,\eta_3(\eta_1))	&=\chi^+(D)\big(2m(\eta_1,F(\eta_1)+\eta_3)+m(F(\eta_1)+\eta_3,F(\eta_1)+\eta_3)\big) \\
& = 2\chi^+(D)m(\eta_1,F(\eta_1))+ \underline{O}(\varepsilon^{3\theta}\nn \eta_1 \nn^3),
\end{align*}
in which we have used the calculations
$$m(\eta_1,\eta_3) =\underline{O}(\varepsilon^{3\theta}\nn \eta_1 \nn^3), \qquad
m(F(\eta_1),\eta_3) = \underline{O}(\varepsilon^{3\theta} \nn \eta_1 \nn^4)$$
(see Proposition \ref{prop:m estimate}, Corollary \ref{cor:F estimates} and Theorem \ref{thm:estimate eta3}) and
$$m(F(\eta_1),F(\eta_1)) = \underline{O}(\varepsilon^{3\theta} \nn \eta_1 \nn^4)$$
(because of \eqref{Formula for m} and Proposition \ref{pm expansion of F}). Observing that
$$
m(\eta_1,F_\mathrm{r}(\eta_1)) = \underline{O}(\varepsilon^{1+2\theta}\nn \eta_1 \nn^3)$$
(see Proposition \ref{prop:m estimate} and Remark \ref{Fr estimate}), we find that
$$\chi^+(D)m(\eta_1,F(\eta_1))=-2\chi^+(D)m(\eta_1^-,(\eta_1^+)^2))+\underline{O}_+(\varepsilon^{3\theta}\nn \eta_1 \nn^3),$$
and it follows from \eqref{Formula for m} and Lemma \ref{lem:approximate identities} that
$$
m(\eta_1^-,(\eta_1^+)^2)
=-\eta_1^-(\eta_1^+)^2 + \underline{O}(\varepsilon^{3\theta} \nn \eta_1 \nn^3).\eqno{\Box}
$$

\begin{proposition} \label{Compute redeq 2}
The estimates
\begin{align*}
\chi^+(D)\KKp_3(\eta_1+F(\eta_1) + \eta_3(\eta_1)) &= -\tfrac{3}{2}\chi^+(D)\big(\eta_1^-(\eta_1^+)^2\big)+\underline{O}_+(\varepsilon^{3\theta}\nn \eta_1 \nn^3), \\
\chi^+(D)\LLp_3(\eta_1+F(\eta_1) + \eta_3(\eta_1)) &= -2\chi^+(D)\big(\eta_1^-(\eta_1^+)^2\big)+\underline{O}_+(\varepsilon^{3\theta}\nn \eta_1 \nn^3)
\end{align*}
hold for each $\eta_1 \in X_1$.
\end{proposition}
\proof
Using the estimates for $F(\eta_1)$ and $\eta_3(\eta_1)$ given in Corollary \ref{cor:F estimates} and Theorem \ref{thm:estimate eta3},
we find that
$$
\KKp_3(\eta_1+F(\eta_1) + \eta_3(\eta_1)) =\KKp_3(\eta_1)+\underline{O}(\varepsilon^{3\theta}\nn \eta_1 \nn^4)
$$
and
$$\chi^+(D)\KKp_3(\eta_1)=-\tfrac{3}{2}\chi^+(D)\big(\eta_1^-(\eta_1^+)^2\big)+\underline{O}_+(\varepsilon^{3\theta}\nn \eta_1 \nn^3))$$
(because of equation \eqref{Formula for KKp3}). It similarly follows from the formula
$$\LLp_3(\eta)= - K_0\eta K_1(\eta)\eta - L_0 \eta L_1(\eta)\eta - \eta_x^2K_0\eta-\eta_x\eta_z L_0\eta + K_2(\eta)\eta$$
and the fact that $K_2(\eta)=m_2(\eta,\eta)$, where $m_2$ is a bounded, symmetric bilinear mapping $\ZZ \times \ZZ \rightarrow \LL(H^{5/2}(\R^2), H^{3/2}(\R^2))$, that
$$
\LLp_3(\eta_1+F(\eta_1) + \eta_3(\eta_1)) =\LLp_3(F(\eta_1)+\eta_1)+\underline{O}(\varepsilon^{3\theta}\nn \eta_1 \nn^4);
$$
using Lemma \ref{KK, LL expansions}(iii) twice yields
$$
\LLp_3(F(\eta_1)+\eta_1)
= \LLp_3(\eta_1) +\underline{O}(\varepsilon^{3\theta} \nn \eta_1 \nn^3)
$$
and
$$
\chi^+(D)\LLp_3(\eta_1)
= -2\chi^+(D)\big(\eta_1^-(\eta_1^+)^2\big) +\underline{O}_+(\varepsilon^{3\theta} \nn \eta_1 \nn^3).\eqno{\Box}
$$

\begin{proposition} \label{Compute redeq 3}
The estimates
\begin{align*}
\KKp_\mathrm{r}(\eta_1+F(\eta_1) + \eta_3(\eta_1)) &= \underline{O}(\varepsilon^{4\theta}\nn \eta_1 \nn^5), \\
\LLp_\mathrm{r}(\eta_1+F(\eta_1) + \eta_3(\eta_1)) &= \underline{O}(\varepsilon^{3\theta}\nn \eta_1 \nn^4)
\end{align*}
hold for each $\eta_1 \in X_1$.
\end{proposition}
\proof
This result follows from Proposition \ref{KK, LL hot estimates}(ii), Corollary \ref{cor:F estimates} and Theorem \ref{thm:estimate eta3}.
\qed\medskip

\begin{proposition} \label{Extra terms}
The estimates
\begin{align*}
-\tfrac{1}{2}\chi^+(D)\big(K_1(\eta_1+F(\eta_1) + \eta_3(\eta_1))(\eta_1+F(\eta_1) + \eta_3(\eta_1))\big)^2 &=\underline{O}_+(\varepsilon^{3\theta}\nn \eta_1 \nn^4), \\
-\tfrac{1}{2}\chi^+(D)\big(L_1(\eta_1+F(\eta_1) + \eta_3(\eta_1))(\eta_1+F(\eta_1) + \eta_3(\eta_1))\big)^2 &=\underline{O}_+(\varepsilon^{3\theta}\nn \eta_1 \nn^4)
\end{align*}
hold for each $\eta_1 \in X_1$.
\end{proposition}
\proof
Using Corollary \ref{cor:F estimates} and Theorem \ref{thm:estimate eta3} we find that
$$-\tfrac{1}{2}\big(K_1(\eta_1+F(\eta_1) + \eta_3(\eta_1))(\eta_1+F(\eta_1) + \eta_3(\eta_1))\big)^2=-\tfrac{1}{2}(K_1(\eta_1)\eta_1)^2 + \underline{O}(\varepsilon^{3\theta}\nn \eta_1 \nn^4),$$
and furthermore
$$-\tfrac{1}{2}\chi^+(D)(K_1(\eta_1)\eta_1)^2 = -\tfrac{1}{2}\chi^+(D)\big((\eta_1\eta_{1x})_x+K_0(\eta_1K_0\eta_1)+L_0(\eta_1 L_0\eta_1)\big)^2=0$$
because of equation \eqref{Formula for K1}.
The second estimate is derived in the same fashion (with equation \eqref{Formula for L1}).
\qed\medskip

\begin{corollary}
The estimate
$$\chi^+(D)\NN_2(\eta_1,\eta_3(\eta_1))=\tfrac{5}{2}\chi^+(D)\big(\eta_1^-(\eta_1^+)^2\big)+\underline{O}_+(\varepsilon^{3\theta}\nn \eta_1 \nn^3) $$
holds for each $\eta_1 \in X_1$.
\end{corollary}

We conclude that the reduced equation for $\eta_1$ is 
$$\eta_1^+-\eta_{1xx}^+-\eta_{1zz}^+-2 K_0\eta_1^+ + 2\varepsilon^2 K_0 \eta_1^+ - \tfrac{11}{2}\chi^+(D)\big(|\eta_1^+|^2\eta_1^+\big) + \underline{O}_+(\varepsilon^{3\theta} \nn \eta_1 \nn^3)=0,$$
which can be further simplified to
$$\eta_1^+-\eta_{1xx}^+-\eta_{1zz}^+-2 K_0\eta_1^+ + 2\varepsilon^2 \eta_1^+ - \tfrac{11}{2}\chi^+(D)\big(|\eta_1^+|^2\eta_1^+\big) + \underline{O}_+(\varepsilon^{3\theta} \nn \eta_1 \nn)=0$$
by an application of Lemma \ref{lem:approximate identities}(iv). Finally, we introduce the nonlinear Schr\"{o}dinger scaling 
$$\eta_1^+(x,z) = \tfrac{1}{2}\varepsilon \zeta(\varepsilon x,\varepsilon z)\ee^{\i x},$$
so that $\zeta \in B_R(0) \subseteq \chi_0(\varepsilon D)H^1(\R^2)$
solves the \emph{perturbed full-dispersion nonlinear Schr\"{o}dinger equation}
\begin{equation}
\varepsilon^{-2}g(e+\varepsilon D)\zeta + 2\zeta - \tfrac{11}{8}\chi_0(\varepsilon D)(|\zeta|^2\zeta) + \varepsilon^{3\theta-2}
\underline{O}_0^\varepsilon(\|\zeta\|_1)=0,
\label{PFDNLS}
\end{equation}
where $R=R_1/\sqrt{2}$ and $e=(1,0)$ (note that $\nn \eta_1 \nn^2 = \|\zeta\|_1^2$ and
the change of variables from $(x,z)$ to $\varepsilon(x,z)$ introduces an additional factor of $\varepsilon$
in the remainder term).
The invariance of the reduced equation under $\eta_1(x,z) \mapsto \eta_1(-x,-z)$ and
$\eta_1(x,z) \mapsto \eta_1(x,-z)$ 
is inherited by \eqref{PFDNLS}, which is invariant under the reflections
$\zeta(x,z) \mapsto \overline{\zeta(-x,z)}$ and $\zeta(x,z) \mapsto \zeta(x,-z)$.

\begin{remark}
In the formal limit $\varepsilon=0$ equation \eqref{PFDNLS} reduces to the nonlinear Schr\"{o}dinger equation
\begin{equation}
-\tfrac{1}{2}\zeta_{xx}-\zeta_{zz} + \zeta -\tfrac{11}{16}|\zeta|^2 \zeta =0. \label{NLS again}
\end{equation}
\end{remark}

\section{Solution of the reduced equation} \label{sec:existence}

In this section we complete our existence theory by proving the following theorem. 

\begin{theorem} \label{Final existence thm}
For each sufficiently small value of $\varepsilon>0$ equation \eqref{PFDNLS} has two small-amplitude
solutions $\zeta^\pm_\varepsilon$ in $\chi_0(\varepsilon D)H^1(\R^2)$ which satisfy
$\zeta^\pm_\varepsilon(x,z) = \overline{\zeta^\pm_\varepsilon(-x,z)}$,
$\zeta^\pm_\varepsilon(x,z) = \zeta^\pm_\varepsilon(x,-z)$ and
$\|\zeta^\pm_\varepsilon - (\pm\zeta_0)\|_1 \lesssim \varepsilon^{1/2}$, where $\zeta_0 \in \SS(\R^2)$ is the unique
symmetric, positive (real) solution of the nonlinear Schr\"{o}dinger equation \eqref{NLS again}.
\end{theorem}

The first step is a result which allows us to `replace' the nonlocal operator in equation \eqref{PFDNLS} with a differential operator.

\begin{proposition}
The inequality
$$\left|\frac{\varepsilon^2}{2\varepsilon^2 + g(e+\varepsilon k)}- \frac{1}{2+k_1^2+2k_3^2}\right| \lesssim \frac{\varepsilon |k|^3}{(1+|k|^2)^2}$$
holds uniformly over $|k| < \delta/\varepsilon$.
\end{proposition}

\proof Clearly
$$\left|\frac{\varepsilon^2}{2\varepsilon^2 + g(e+\varepsilon k)}- \frac{1}{2+k_1^2+2k_3^2}\right|
=
\frac{|g(e+\varepsilon k)-\varepsilon^2(k_1^2+2k_3^2)|}{(2\varepsilon^2 + g(e+\varepsilon k))(2+k_1^2+2k_3^2)},
$$
while
$$g(e+s)-s_1^2-2s_2^2 \lesssim |s|^3, \qquad |s| \leq \delta$$
and
$$g(e+s) \gtrsim |s|^2, \qquad s \in \R^2.$$
It follows that
$$\left|\frac{\varepsilon^2}{2\varepsilon^2 + g(e+\varepsilon k)}- \frac{1}{2+k_1^2+2k_3^2}\right|
\lesssim
\frac{\varepsilon |k|^3}{(1+|k|^2)^2}, \qquad |k| < \delta/\varepsilon.\eqno{\Box}$$

Using this proposition, one can write equation \eqref{PFDNLS} as
\begin{equation}
\zeta+F_\varepsilon(\zeta)=0, \label{2nd zeta eqn}
\end{equation}
where
$$
F_\varepsilon(\zeta)=-\tfrac{11}{16}\left(1-\tfrac{1}{2}\partial_x^2 - \partial_z^2\right)^{-1}\chi_0(\varepsilon D)\big(|\zeta|^2\zeta\big)
+ \varepsilon^{1/2}\underline{O}_1^\varepsilon(\|\zeta\|_1)
$$
and we have chosen the concrete value $\theta =5/6$, so that $\varepsilon^{3\theta-2}=\varepsilon^{1/2}$.
It is convenient to replace equation \eqref{2nd zeta eqn} with
\begin{equation}
\zeta+\tilde{F}_\varepsilon(\zeta)=0, \label{3rd zeta eqn}
\end{equation}
where $\tilde{F}_\varepsilon(\zeta)=F_\varepsilon(\chi_0(\varepsilon D)\zeta)$ and study it in the fixed space
$H^1(\R^2)$ (the solution sets of \eqref{2nd zeta eqn} and \eqref{3rd zeta eqn} evidently coincide).
Equation \eqref{3rd zeta eqn} is solved using the following version of the implicit-function theorem.

\begin{theorem} \label{IFT}
Let $\XX$ be a Banach space, $X_0$ and $\Lambda_0$ be open neighbourhoods of respectively $x^\star$ in $\XX$ and the origin in $\R^n$
and $G:  X_0 \times \Lambda_0 \rightarrow \XX$ be a function which is differentiable with respect to $x \in X_0$ for each $\lambda \in \Lambda_0$.
Suppose that 
$G(x^\star,0)=0$, $\mathrm{d}_1G[x^\star,0]: \XX \rightarrow \XX$ is an isomorphism,
$$\lim_{x \rightarrow x^\star}\|\mathrm{d}_1G[x, 0]-\mathrm{d}_1G[x^\star,0]\|_{\LL(\XX)}=0$$
and
$$\lim_{\lambda \rightarrow 0} \|G(x,\lambda)-G(x,0)\|_{\XX}=0, \quad \lim_{\lambda \rightarrow 0} \
\|\mathrm{d}_1G[x,\lambda]-\mathrm{d}_1G[x,0]\|_{\LL(\XX)}=0$$
uniformly over $x \in X_0$.

There exist open neighbourhoods $X$ of $x^\star$ in $\XX$ and $\Lambda$ of $0$ in $\R^n$
(with $X \subseteq X_0$, $\Lambda \subseteq \Lambda_0)$ and a uniquely determined mapping
$h: \Lambda \rightarrow X$ with the properties that
\begin{itemize}
\item[(i)]
$h$ is continuous at the origin (with $h(0)=x^\star$),
\item[(ii)]
$G(h(\lambda),\lambda)=0$ for all $\lambda \in \Lambda$,
\item[(iii)]
$x=h(\lambda)$ whenever $(x,\lambda) \in X \times \Lambda$ satisfies $G(x,\lambda)=0$.
\end{itemize}

Furthermore, the existence of $\alpha>0$ such that $\|G(x,\lambda)-G(x,0)\|_{\XX} \lesssim |\lambda|^\alpha$ for all $\lambda \in \Lambda_0$ and $x \in X_0$ implies that $\|h(\lambda)-h(0)\|_{\XX} \lesssim |\lambda|^\alpha$ for all $\lambda \in \Lambda$.
\end{theorem}

We establish Theorem \ref{Final existence thm} by applying Theorem \ref{IFT} with
$$\XX=H_\mathrm{e}^1(\R^2,\C)=\{\zeta \in H^1(\R^2,\C): \zeta(x,z) = \overline{\zeta(-x,z)},
\ \zeta(x,z) = \zeta(x,-z)\},$$
$X=B_R(0)$, where $R$ is chosen large enough that $\zeta_0 \in X$,
$\Lambda_0=(-\varepsilon_0,\varepsilon_0)$ for a sufficiently small value of $\varepsilon_0$
and
$$G(\zeta,\varepsilon):=\zeta+\tilde{F}_{|\varepsilon|}(\zeta)$$
(here $\varepsilon$ is replaced by $|\varepsilon|$ so that $G(\zeta,\varepsilon)$ is defined for $\varepsilon$ in a full neighbourhood of
the origin in $\R$).

Observe that
\begin{align*}
G(\zeta&,\varepsilon)-G(\zeta,0) \\
&= -\tfrac{11}{16}\left(1-\tfrac{1}{2}\partial_x^2 - \partial_z^2\right)^{-1}\bigg(\chi_0(|\varepsilon| D)\Big(|\chi_0({|\varepsilon|} D)\zeta|^2\chi_0({|\varepsilon|} D)\zeta\Big)-|\zeta|^2\zeta\bigg)
 + {|\varepsilon|}^{1/2}\underline{O}_1^{|\varepsilon|}(\|\zeta\|_1) \\
 &= -\tfrac{11}{16}\left(1-\tfrac{1}{2}\partial_x^2 - \partial_z^2\right)^{-1}\bigg(\chi_0(|\varepsilon| D)\Big(
|\chi_0({|\varepsilon|}D)\zeta|^2\big(\chi_0({|\varepsilon|} D)-I\big)\zeta +|\zeta|^2\big(\chi_0({|\varepsilon|} D)-I\big)\zeta\\
&\hspace{2.7in}\mbox{}+
\zeta\chi_0({|\varepsilon|}D)\zeta\big(\chi_0({|\varepsilon|} D)-I\big)\overline{\zeta} \Big)\\
 & \hspace{2.1in} \mbox{}+\Big(\chi_0({|\varepsilon|} D)-I\Big)|\zeta|^2\zeta\bigg)
 + {|\varepsilon|}^{1/2}\underline{O}_1^{|\varepsilon|}(\|\zeta\|_1).
 \end{align*}
Noting that
$$\|\chi_0(|\varepsilon| D) - I\|_{{\mathcal L}(H^1({\mathbb R}^2,\C), H^{1/2}({\mathbb R}^2,\C))} \lesssim |\varepsilon|^{1/2}$$
because
\begin{align*}
\| \chi_0(|\varepsilon| D) u - u \|_{1/2}^2 &= \int_{|k| > \frac{\delta}{|\varepsilon|}} (1+|k|^2)^{1/2} |\hat{u}|^2 \dk \\
& \leq \sup_{|k| > \frac{\delta}{|\varepsilon|}} (1+|k|^2)^{-1/2} \int_{|k| > \frac{\delta}{|\varepsilon|}} (1+|k|^2) |\hat{u}|^2 \dk\\
& \leq \frac{1}{\big(1+\frac{\delta^2}{|\varepsilon|^2}\big)^{1/2}} \|u\|_1^2,
\end{align*}
and similarly
$$\|\chi_0(|\varepsilon| D) - I\|_{{\mathcal L}(H^{1/2}({\mathbb R}^2,\C), L^2({\mathbb R}^2,\C))} \lesssim |\varepsilon|^{1/2},$$
and that pointwise multiplication defines bounded trilinear mappings $H^1(\R^2,\C)^3 \rightarrow
H^{1/2}(\R^2,\C)$ and $H^1(\R^2,\C)^2 \times H^{1/2}(\R^2,\C) \rightarrow L^2(\R^2,\C)$
(see H\"{o}rmander \cite[Theorem 8.3.1]{Hoermander}), we find that
$$\|G(\zeta,\varepsilon)-G(\zeta,0)\|_1 \lesssim |\varepsilon|^{1/2}$$
uniformly over $\zeta \in B_R(0)$. Here we have also used the estimate $\|\chi_0(|\varepsilon| D)u \|_s \leq \|u\|_s$
for all $u \in H^s(\R^2,\C)$ and
and the fact that $\left(1-\tfrac{1}{2}\partial_x^2 - \partial_z^2\right)^{-1}$ maps
$L^2({\mathbb R}^2,\C)$ continuously into $H^1({\mathbb R}^2,\C)$. A similar calculation shows that
$$\|\mathrm{d}_1G[\zeta,\varepsilon]-\mathrm{d}_1G[\zeta,0]\|_{\LL(H^1(\R^2,\C))} \lesssim |\varepsilon|^{1/2}$$
uniformly over $\zeta \in B_R(0)$.

Furthermore the equation
\begin{equation}
G(\zeta,0)=\zeta-\tfrac{11}{16}\left(1-\tfrac{1}{2}\partial_x^2 - \partial_z^2\right)^{-1}|\zeta|^2\zeta=0 \label{Limiting eqn}
\end{equation}
has a unique symmetric, positive (real) solution $\zeta_0 \in \SS(\R^2,\C)$ (see Sulem \& Sulem \cite[\S 4.2]{SulemSulem} and the references therein).
The fact that $\mathrm{d}_1G[\pm\zeta_0,0]$ is an isomorphism is conveniently established by using real
coordinates. Define
$\zeta_1=\re \zeta$ and $\zeta_2=\im \zeta$, so that
$$\mathrm{d}_1G[\pm\zeta_0,0](\zeta_1+\mathrm{i}\zeta_2)=G_1(\zeta_1) + \mathrm{i} G_2(\zeta_2),$$
where $G_1: H_\mathrm{e}^1(\R^2,\R) \rightarrow H_\mathrm{e}^1(\R^2,\R)$ and
$G_2: H_\mathrm{o}^1(\R^2,\R) \rightarrow H_\mathrm{o}^1(\R^2,\R)$ are defined by
$$G_1(\zeta_1)=\zeta_1 - \tfrac{33}{16}\left(1-\tfrac{1}{2}\partial_x^2 - \partial_z^2\right)^{-1} \zeta_0^2 \zeta_1, \qquad
G_2(\zeta_2)=\zeta_2 - \tfrac{11}{16}\left(1-\tfrac{1}{2}\partial_x^2 - \partial_z^2\right)^{-1} \zeta_0^2 \zeta_2$$
with
\begin{align*}
H^n_\mathrm{e}(\R^2,\R) &= \{\zeta_1 \in H^n(\R^2,\R): \zeta_1(x,z)=\zeta_1(-x,z),\ \zeta_1(x,z)=\zeta_1(x,-z)\}, \\
H^n_\mathrm{o}(\R^2,\R) &= \{\zeta_2 \in H^n(\R^2,\R): \zeta_2(x,z)=-\zeta_2(-x,z),\ \zeta_2(x,z)=\zeta_1(x,-z)\}
\end{align*}
for $n \in {\mathbb N}_0$.
The formulae
$$\zeta_1 \mapsto \tfrac{33}{16}\left(1-\tfrac{1}{2}\partial_x^2 - \partial_z^2\right)^{-1} \zeta_0^2 \zeta_1, \qquad
\zeta_2 \mapsto \tfrac{11}{16}\left(1-\tfrac{1}{2}\partial_x^2 - \partial_z^2\right)^{-1} \zeta_0^2 \zeta_2$$
define compact operators $H^1(\R^2,\R) \rightarrow H^1(\R^2,\R)$,
$H^1_\mathrm{e}(\R^2,\R) \rightarrow H^1_\mathrm{e}(\R^2,\R)$ and $H^1_\mathrm{o}(\R^2,\R) \rightarrow H^1_\mathrm{o}(\R^2,\R)$,
so that $G_1$, $G_2$ are Fredholm operators with index $0$. Writing
$$T_1\zeta_1 = \zeta_1 - \tfrac{1}{2}\zeta_{1xx}-\zeta_{1zz} - \tfrac{33}{16}\zeta_0^2\zeta_1, \qquad
T_2\zeta_2 =  \zeta_2 - \tfrac{1}{2}\zeta_{2xx}-\zeta_{2zz} - \tfrac{11}{16}\zeta_0^2 \zeta_2,$$
we find that the kernels of $G_1$ and $G_2$ coincide with respectively the kernels of the linear operators
$T_1: H^2_\mathrm{e}(\R^2,\R) \subseteq L^2_\mathrm{e}(\R^2,\R) \rightarrow L^2_\mathrm{e}(\R^2,\R)$
and
$T_2: H^2_\mathrm{o}(\R^2,\R) \subseteq L^2_\mathrm{o}(\R^2,\R) \rightarrow L^2_\mathrm{o}(\R^2,\R)$.
It is however known that the kernels of $T_1$, $T_2: H^2(\R^2,\R) \subseteq L^2(\R^2,\R) \rightarrow L^2(\R^2,\R)$
are respectively $\langle\zeta_{0x}, \zeta_{0z}\rangle$ and $\langle\zeta_0\rangle$  (see Chang \emph{et al.}\ 
\cite{ChangGustafsonNakanishiTasi07}). The kernels of $G_1$, $G_2$ are therefore trivial,
so that $G_1$, $G_2$ and hence $\mathrm{d}_1G[\pm\zeta_0,0]$ are isomorphisms.

It remains to confirm that tracing back the changes of variable
$$\eta = \eta_1 + F(\eta_1) + \eta_3(\eta_1),
\quad \eta_1=\eta_1^++\overline{\eta_1^+},
\quad \eta_1^+(x,z) = \tfrac{1}{2}\zeta^\pm_\varepsilon(\varepsilon x, \varepsilon z)\ee^{\ii x}$$
leads to the estimate
$$\eta(x,z)=\pm \varepsilon \zeta_0(\varepsilon x,\varepsilon z)\cos x + o(\varepsilon)$$
uniformly over $(x,z) \in \R^2$. The key is to show that
$$\|\zeta^+_\varepsilon-\zeta_0\|_\infty \lesssim \varepsilon^\Delta$$
for any $\Delta \in (0,1/2)$; here we choose the concrete value $\Delta=1/4$. This result follows from the calculation
\begin{align*}
\|\zeta_\varepsilon^+-\zeta_0\|_\infty &\lesssim \|\zeta_\varepsilon^+-\zeta_0\|_{5/4} \\
& = \|(1+|k|^2)^{5/8}(\hat{\zeta}_\varepsilon^+-\hat{\zeta}_0)\|_{L^2(|k|<\delta/\varepsilon)}+
\|(1+|k|^2)^{5/8}\hat{\zeta}_0\|_{L^2(|k|>\delta/\varepsilon)}
\end{align*}
(because the support of $\hat{\zeta}_\varepsilon$ lies in $\overline{B}_{\delta/\varepsilon}(0)$) and
\begin{align*}
\|(1+|k|^2)^{5/8}(\hat{\zeta}_\varepsilon^+-\hat{\zeta}_0)\|_{L^2(|k|<\delta/\varepsilon)}
& \lesssim
\varepsilon^{-1/4}\|(1+|k|^2)^{1/2}(\hat{\zeta}_\varepsilon^+-\hat{\zeta}_0)\|_{L^2(|k|<\delta/\varepsilon)} \\
& \leq
\varepsilon^{-1/4}\|(1+|k|^2)^{1/2}(\hat{\zeta}_\varepsilon^+-\hat{\zeta}_0)\|_0 \\
& =
\varepsilon^{-1/4}\|\zeta_\varepsilon^+-\zeta_0\|_1, \\
& \lesssim \varepsilon^{1/4}, \\
\\
\|(1+|k|^2)^{5/8}\hat{\zeta}_0^+\|_{L^2(|k|>\delta/\varepsilon)}^2
& = \int_{|k|>\frac{\delta}{\epsilon}} (1+|k|^2)^{5/4} |\hat{\zeta}_0|^2 \\
& \lesssim \varepsilon
\end{align*}
(because $\hat{\zeta}_0 \in \SS({\mathbb R}^2)$, so that in particular $|\hat{\zeta}_0(|k|)|^2 \lesssim (1+|k|^2)^{-11/4}$).
It follows that
\begin{align*}
\eta_1^+(x,z) & = \varepsilon \zeta_0(\varepsilon x,\varepsilon z)\ee^{\ii x}
+\tfrac{1}{2}\varepsilon (\zeta_\varepsilon^+-\zeta_0)(\varepsilon x,\varepsilon z)\ee^{\ii x}  \\
& = \varepsilon \zeta_0(\varepsilon x,\varepsilon z)\ee^{\ii x} + O(\varepsilon^{5/4}),
\end{align*}
uniformly in $(x,z)$. (These estimates remain valid when $\zeta_\varepsilon^+$ and
$\zeta_0$ are replaced by respectively $\zeta_\varepsilon^-$ and $-\zeta_0$.)

Furthermore
$$\|\eta_3(\eta_1)\|_\infty \lesssim \|\eta_3(\eta_1)\|_3 \lesssim \varepsilon^{10/6} \nn \eta_1 \nn^2 \lesssim \varepsilon^{10/6}$$
by Theorem \ref{thm:estimate eta3} (recall that we have chosen $\theta=5/6)$, while
$$\|F(\eta_1)\|_\infty=O(\varepsilon^{11/6})$$
because
$$F(\eta_1) = -2\left((\eta_1^+)^2+\big(\overline{\eta_1^+}\big)^2\right)+F_\mathrm{r}(\eta_1),$$
where
$$\|F_\mathrm{r}(\eta_1)\|_\infty \lesssim \|F_\mathrm{r}(\eta_1)\|_3 \lesssim \|F_\mathrm{r}(\eta_1)\|_1 \lesssim \varepsilon^{11/6}\nn \eta_1\nn^2 \lesssim \varepsilon^{11/6}$$
(see Proposition \ref{pm expansion of F}; the second estimate follows by the fact that the support of $\FF[F_\mathrm{r}(\eta_1)]$
is bounded independently of $\varepsilon$).

\noindent\\
{\bf Acknowledgement.} E. Wahl\'{e}n was supported by the Swedish Research Council, grant no.\ 2016-04999.

\end{document}